\documentclass[11pt]{article}
\usepackage[english]{babel}
\usepackage[latin1]{inputenc}
\usepackage{amsmath,amssymb,enumerate,color}
\usepackage{ulem}
\usepackage{psfig,labelfig}
\usepackage{graphicx}
\usepackage[dvips]{epsfig}
\usepackage[T1]{fontenc}

\newtheorem{thm}{Theorem}[section]
\newtheorem{cor}[thm]{Corollary}
\newtheorem{lem}[thm]{Lemma}
\newtheorem{cla}[thm]{Claim}

\newtheorem{defi} [thm]{Definition}

\newcommand{\R}{{\mathbb R}}

\DeclareMathOperator{\arcsinh}{arcsinh}
\DeclareMathOperator{\arccosh}{arccosh}
\DeclareMathOperator{\arctanh}{arctanh}

\DeclareMathOperator{\dist}{dist}
\DeclareMathOperator{\capa}{cap}
\DeclareMathOperator{\area}{area}
\DeclareMathOperator{\ao}{\alpha_1}
\DeclareMathOperator{\at}{\alpha_2}

\title{On the second successive minimum of the Jacobian of a Riemann surface}

\author{Muetzel, Bjoern\thanks{E-mail address : bjorn.mutzel@epfl.ch}\\
 \\
\small Department of Mathematics, Ecole Polytechnique F\'{e}d\'{e}rale de Lausanne, Station 8, \\
\small CH-1015 Lausanne, Switzerland \\[-0.8ex]
 }

\begin{document}
\maketitle

\begin{abstract}
To a compact Riemann surface of genus $g$ can be assigned a \textit{principally polarized abelian variety} (PPAV) of dimension $g$, the \textit{Jacobian} of the Riemann surface. The Schottky problem is to discern the Jacobians among the PPAVs. Buser and Sarnak showed, that the square of the first successive minimum, the squared norm of the shortest non-zero vector in the lattice of a Jacobian of a Riemann surface of genus $g$ is bounded from above by $\log(4g)$, whereas it can be of order $g$ for the lattice of a PPAV of dimension $g$. We show, that in the case of a hyperelliptic surface this geometric invariant is bounded from above by a constant and that for any surface of genus $g$ the square of the second successive minimum is equally of order $\log(g)$. We obtain improved bounds for the $k$-th successive minimum of the Jacobian, if the surface contains small simple closed geodesics.
\end{abstract}

\section{Introduction}

 A \textit{principally polarized abelian variety (PPAV)} of dimension $g$ may be defined as a pair $(A,H)$, where $A={\mathbb{C}}^g / L$ is a complex torus of dimension $g$, the quotient of $ {\Bbb C}^g $ modulo a lattice $L$. Furthermore, $H$, the \textit{polarization}, is a positive definite hermitian form, whose imaginary part $ImH$ is integral on the lattice points of the lattice $L$. It is \textit{principal}, if it satisfies certain conditions (see \cite{bl}, Section 3.1.).\\
Certain PPAV arise from compact Riemann surfaces in the following way. Let $\left( {\alpha _i } \right)_{i = 1...2g}$ be a \textit{canonical homology basis} of closed cycles on the surface $S$. This basis is assumed to be given in the following way. Each $\alpha _i$ is a simple closed curve and the curves are paired, such that for each $\alpha _i$ there exists exactly one $\alpha _{\tau(i)} \in \left( {\alpha _j } \right)_{j = 1...2g}$, that intersects $\alpha_i$ in exactly one point and there are no other intersection points. In the vector space of harmonic forms on $S$ let $\left( {u _k } \right)_{k = 1...2g}$ be a \textit{dual basis for} $\left( {\alpha _i } \right)_{i = 1...2g}$ defined by

\[
\int\limits_{\alpha_i} {u_k } = \delta_{ik} .
\]

The \textit{period matrix} $P_S$ of a compact Riemann surface (R.S.) $S$  of genus $g$ is the Gram matrix

\begin{equation}
P_S=\left(\left\langle {u_i,u_j} \right\rangle  \right)_{i,j= 1...2g}=\left( {\int\limits_S {u_i  \wedge {}^ * } u_j } \right)_{i,j= 1...2g}   \nonumber
\end{equation}

This period matrix $P_S$ defines a complex torus, the \textit{Jacobian} or \textit{Jacobian variety} $J(S)$ of the Riemann surface $S$. By Riemann's period relations the Jacobian is a PPAV.\\
The moduli or parameter space $\mathcal{A}_g$ for PPAVs of dimension $g$ has dimension $\tfrac{1}{2}g(g + 1)$, whereas the moduli space of compact Riemann surfaces of genus $g$, $\mathcal{M}_g$, has dimension $3g-3$. The assignment of the Jacobian $J(S)$ to the R.S. $S$ provides a mapping $t: \mathcal{M}_g \to \mathcal{A}_g$. By Torelli's theorem, this mapping is injective. In general, it is not known, if a given PPAV is the image of a Jacobian under $t$. The Schottky problem is to describe the sublocus $t(\mathcal{M}_g)$ in $\mathcal{A}_g$. \\
The closure of the sublocus of Jacobian varieties $\overline{t(\mathcal{M}_g)}$ in the parameter space $\mathcal{A}_g$ is only equal to $\mathcal{A}_g$ for $g=2$ and $3$. For $g \geq 4$ it is a proper closed subset. Several analytic approaches have been used to further characterize $t(\mathcal{M}_g)$. Notably van Gemen proved in \cite{vg}, that  $t(\mathcal{M}_g)$ is an irreducible component of the locus $\mathcal{S}_g$ defined by the Schottky-Jung polynomials. However, an exact description of the locus $t(\mathcal{M}_g)$, given in terms of polynomials of theta constants that vanish on $t(\mathcal{M}_g)$ but not on $\mathcal{A}_g$, is only known for $g=4$ (\cite{sc}). Shiota \cite{sh} showed, that an indecomposable PPAV is the Jacobian of a Riemann surface, if the corresponding theta function fulfills the K-P differential equation. However, a solution to this equation can not as yet be determined explicitly.\\
In \cite{bs}, it was shown, that the Jacobians can be characterized among the PPAVs with the help of a geometric invariant of the lattice of the PPAVs, the first successive minimum or shortest non-zero lattice vector, whose square is also called the \textit{minimal period length} of the PPAV. Here the \textit{$k$-th successive minimum of a complex lattice $L$} of dimension $g$ is defined by
\[
m_k(L)=\min \left\{ {r \in \mathbb{R}^ +  \left| {\exists \left\{ {l_1 ,...,l_k } \right\} \subset L,\text{ lin. independent over } \mathbb{R} , \left\| {l_i } \right\| \leqslant r\forall i } \right.} \right\}
\]

The \textit{k-th successive minimum of a PPAV} $({A=\Bbb C}^g/L,H)$, $m_k(A,H)$, is defined as the $k$-th successive minimum of its lattice $L$. Here the norm $\left\| {\cdot } \right\|=\left\| {\cdot } \right\|_H$ is the norm induced by the hermitian form $H$. If $\left(l_i \right)_{i \in \{1..2g\}}$ is a lattice basis of $L$, then the corresponding Gram matrix
\[
P_H=\left({\left\langle {l_i,l _j} \right\rangle }_H \right)_{i,j= 1...2g}
\]

has determinant $1$, due to the fact, that the PAV is principal. Therefore we can apply the general upper bounds on the successive minima from Minkowski's theorems (see \cite{gl}) to the case of a PPAV $(A,H)$ of dimension $g$, whereas the lower bound for Hermite's constant over the PPAVs \\

$$\delta_{2g}=\mathop{max} \limits_{(A,H) \in \mathcal{A}_g} m_1(A,H)^2 .$$

was proven in \cite{bs} :

\[
\frac{g}
{{\pi e}} \approx \frac{1}
{\pi }\sqrt[g]{{2g!}}{\text{  }} \leqslant \delta_{2g}  \leqslant  {\text{  }}\frac{4}
{\pi }\sqrt[g]{{g!}} \approx \frac{{4g}}
{{\pi e}}
\]

The approximation applies to large $g$. By Minkowski's second theorem we have

\[
\prod\limits_{{\text{k = 1}}}^{{\text{2g}}} {{\text{m}}_{\text{k}} {\text{(A,H)}}^{\text{2}} } {\text{  }}\mathop  \leqslant \limits^{} {\text{  }}\left( {\frac{4}
{\pi }} \right)^g g!^2
\]

Buser and Sarnak showed in \cite{bs}, that the shortest non-zero lattice vector of the Jacobian of a compact Riemann surface of genus $g$ can be maximally of order $\log(g)$  :

\begin{thm}
If $\eta_{2g}=\max \limits_{(A,H) \in t(\mathcal{M}_g)} m_1(A,H)^2 $ , then
\[
c \log g \leq \eta_{2g} \leq \frac{3}{\pi} \log(4g-2),
\]
where $c$ is a constant.
\label{thm:bs}
\end{thm}

To extend this theorem, we are going to prove the following :

\begin{thm}
Let $S$ be a compact R.S. of genus $g$ and let $J(S)$ be its Jacobian. Then
\[
m_1(J(S))^2  \leq \log(4g-2)  \text{ \ \ and \ \ }  m_2(J(S))^2 \leq 3.1 \log(8g-7)
\]
\label{thm:main}
\end{thm}
For the second successive minimum of a PPAV $(A,H)$ of dimension $g$ we obtain by Minkowski's second theorem :
\[
m_2 (A,H)^2  \leqslant \frac{1}
{\sqrt[g]{m_1 (A,H)}}\left( {\frac{{4\sqrt[g]{{g!}}}}
{\pi }} \right)^{2g/(2g - 1)}  \approx    \frac{{4g}}
{\sqrt[g]{m_1 (A,H)}{\pi e}},
\]

where the approximation applies for large $g$. Furthermore there exist examples of PPAVs where $m_2(A,H)^2$ is of order $g$. This follows from the fact, that PPAVs, whose shortest lattice vector is maximal, have a basis of minimal non-zero vectors (see \cite{ber}). In this case all $m_k(A,H)^2$ are of order $g$. In contrast, we have for the Jacobian of a Riemann surface, $J(S)$, that $m_1(J(S))^2$ and $m_2(J(S))^2$ are both of order $\log(g)$, independent of the length of the shortest non-zero lattice vector.\\
If a R.S. contains a certain number of mutually disjoint small simple closed geodesics, we obtain the following corollary of \textbf{Theorem~\ref{thm:main}} :

\begin{cor}
Let $S$ be a compact R.S. of genus $g$, that contains $n$ disjoint simple closed geodesics $\left(\eta_j\right)_{j=1,..,n}$ of length smaller than $t$. If we cut open $S$ along these geodesics, then the decomposition contains $m$ R.S. $S_i$ of signature $(g_i,n_i)$, with $g_i > 0$. There exist $m$ linear independent vectors $\left(l_i\right)_{i=1,..,m}$ in the lattice of the Jacobian $J(S)$, such that
\[
{\left\| {l_i } \right\|_H}^2 \leq \frac{(n_i+1)\max \{ 4\log(4g_i+2n_i-3),t\}}{\pi - 2\arcsin(M)} \text{ \ \ for \ \ } i \in \{1,...,m\},
\]
where $M=\min \{ \frac{\sinh(\frac{t}{2})}{\sqrt{\sinh(\frac{t}{2})^2 + 1}},\frac{1}{2}\}$
\label{thm:cor_main}
\end{cor}

As the vectors $l_i$ are linearly independent, the corollary implies improved bounds for a certain number of $m_k(J(S))$. This corollary is related to a Theorem of Fay. In \cite{fa}, chap. III a sequence of Riemann surfaces $S_t$ is constructed, where $t$ denotes the length of a separating simple closed geodesic $\eta$. $\eta$ divides $S_t$ into two surfaces $S_i$ of signature $(g_i,1), i \in \{1,2\}$. If $t \rightarrow 0$ then the period matrix for a suitable canonical homology basis converges to a block matrix, where each block is in $M_{2g_i}(\R)$.\\
If $\eta$ is any separating geodesic, that separates a R.S. $S$ into two surfaces $S_i$ of signature $(g_i,1), i \in \{1,2\}$.is small enough. Applying \textbf{Lemma~\ref{thm:lem_bound2}}, we obtain a slightly better bound than in the corollary :
\[
m_i(J(S))^2 \leq \log(8g_i-2)  \text{ \ \ for \ \ } i \in \{1,2\}.
\]
The corollary shows, that indeed the first two successive minima of the Jacobian of the surfaces $S_i$ can only be of order $log(g_i)$ and gives explicit bounds depending on the length of $t$.\\
It shows, that $m_1(J(S))^2$ and $m_2(J(S))^2$ of a R.S. with a sufficiently small separating scg, is at most of the order of the first successive minimum of a R.S. of genus $g_1$ and $g_2$, respectively.

Using the same methods as in \cite{bs}, we furthermore show, that

\begin{thm}
If $S$ is a hyperelliptic  R.S. of genus $g$ and $J(S)$ its Jacobian, then
\[
    m_1(J(S))^2 \leq \frac{3\log(3+2\sqrt{3}+2\sqrt{5+3\sqrt{3}})}{\pi} = 2.4382...
\]
\label{thm:hyper}
\end{thm}

It is worth mentioning, that this  result follows from a simple refinement of the proof, that the systole of hyperelliptic surfaces is bounded from above by a constant, which was shown in \cite{ba} and \cite{je}. \\

\section{Relating the length of lattice vectors of the Jacobian to geometric data of the surface}

In \cite{bs} an upper bound for the norm of a certain lattice vector of a Jacobian of a surface $S$ is obtained by linking the norm of the vector to the length of a non-separating simple closed geodesic on $S$ and the width of its collar, a topological tube around this geodesic. This approach can be further expanded.\\
For simplification the following expressions will be abbreviated. A simple closed geodesic will be denoted \textit{scg} and a non-separating simple closed geodesic \textit{nsscg} and a separating simple closed geodesic \textit{sscg}. By abuse of notation we will denote the length of a geodesic arc by the same letter as the arc itself, if it is clear from the context.\\
Let $S$ be a compact R.S. and $\left(\alpha_i\right)_{i= 1...2g}$ a canonical homology basis on $S$. The collar of a scg $\gamma$, $C(\gamma)$, is defined by

\[
C(\gamma)=\left\{x \in S \mid \dist(x,\gamma) < w \right\}.
\]

Here $w$ is the supremum of all $\omega$, such that the geodesic arcs of length $\omega$ emanating  perpendicularly from  $\gamma$ are pairwise disjoint. For a given $\alpha_i$, let $\alpha_{\tau(i)}$ be the unique scg in the canonical homology basis, that intersects $\alpha_i$. As in \cite{bs}, p. 36 test forms $u'_i$ may be defined on the collar of an $\alpha_i$, that satisfy

\begin{equation}
\int\limits_{\alpha_{j}} {u'_{i} } = \left\{ \begin{array}{ll}
1 & \mbox{if $j=\tau(i)$}\\
0 & \mbox{if $j \neq \tau(i)$}.\\
\end{array} \right.
\label{eq:uprime}
\end{equation}

Among all differential forms on $S$, that satisfy equation (\ref{eq:uprime}), the corresponding harmonic form $u_{\tau(i)}$ in the dual basis for the homology basis $\left(\alpha_i\right)_{i= 1...2g}$ is minimizing with respect to the scalar product $\left\langle {\cdot,\cdot} \right\rangle= \int\limits_S {\cdot  \wedge {}^ * \cdot }$. Therefore

\begin{equation}
      \int\limits_S {u_{\tau(i)}  \wedge {}^ * u_{\tau(i)} } \leq  \int\limits_S {u'_{i}  \wedge {}^ * u'_{i} } \text{ \ \ for all \ \ } i \in \left\{ 1...2g \right\} .
\label{eq:utau}
\end{equation}

If $(A,H)$, where $A={\mathbb{C}}^g / L$ is the Jacobian $J(S)$ of the surface $S$, then

\[
\left({\left\langle {l_i,l _j} \right\rangle }_H \right)_{i,j= 1...2g} = P_H = P_S =\left( {\int\limits_S {u_i  \wedge {}^ * } u_j } \right)_{i,j= 1...2g}
\]

by Riemann's period relations and therefore we have for all $i$

\[
{\left\langle {l_i,l _i} \right\rangle }_H  =  \int\limits_S {u_i  \wedge {}^ * u_{i} }.
\]

The $\left(l_i\right)_{i= 1...2g}$ are linear independent vectors of the lattice $L$ and if we can obtain an upper bound for the test forms $\left(u'_i\right)_{i= 1...k}$, we obtain an upper bound on $m_k(A,H)^2$ by equation (\ref{eq:utau}),

\begin{equation}
m_k(A,H)^2 \leq \max \limits_{i \in \left\{ {1,...,k} \right\}}  \int\limits_S {u'_i  \wedge {}^ * u'_{i} }
\label{eq:mk1}
\end{equation}

The capacity of $C(\alpha_i)$ (see \cite{bs}, p. 36), $\capa (C(\alpha_i))$ provides this upper bound for the squared norm of an $u'_i$. The bound is given by
\begin{equation}
\left\langle {u'_i,u'_i} \right\rangle \leq \capa (C(\alpha_i)) = \frac{l(\alpha_i)}{\pi-2 \cdot \arcsin(\frac{1}{\cosh(w_i)})},
\label{eq:mk2}
\end{equation}

where $w_i$ denotes the width of the collar $C(\alpha_i)$.  $\capa (C(\alpha_i))$ is a strictly increasing function with respect to $w_i$. The following values $W$ and $W'$ for the width of a collar occur frequently in our proof :
\[
    W=\arccosh(2)=1.3169...    \text{ \ \ and \ \ } W'= \arctanh(2/3)=0.8047...
\]
If $w_i=W$, we have that
\[
    \capa (C(\alpha_i))= \frac{3 l(\alpha_i)}{2\pi}  \leq 0.5 l(\alpha_i).
\]
If $w_i=W'$, we obtain that
\[
    \capa (C(\alpha_i))= \frac{l(\alpha_i)}{\pi-2\arcsin(\frac{\sqrt{5}}{3})}  \leq 0.7 l(\alpha_i).
\]
The upper bound in \textbf{Theorem~\ref{thm:bs}} follows from the fact, that a canonical homology basis $\left(\alpha_i\right)_{i= 1...2g}$ can always be constructed, such that $\ao$ is the shortest nsscg on a Riemann surface $S$. It was shown in \cite{bs}, that the length of the shortest nsscg $\ao$ is smaller than $2 \log(4g-2)$ for any R.S. of genus $g$ and that its collar width $w_1$ is bounded from below by $W'$. It follows from the above equations, that $m_1(J(S))^2$ is bounded from above. A more refined analysis shows, that $m_1(J(S))^2 \leq \frac{3}{\pi} \log(4g-2)$.\\
We obtain \textbf{Theorem~\ref{thm:main}} by showing, that there exist two short nsscg, $\ao$ and $\at$, whose collar widths are bounded from below and which can be incorporated together into the canonical homology basis. In principle this approach would provide further bounds for the consecutive $m_k(J(S))$, but finding bounds for both collar width and length of the nsscgs has already proven to be very technical for $k=2$.\\
\textbf{Theorem~\ref{thm:hyper}} follows from the fact, that the non-separating systole of a hyperelliptic surface is bounded from above by a constant, independent of the genus.

\section{General upper bounds for the length of short scgs on a Riemann surface}

\begin{figure}[h!]
\SetLabels
\L(.14*.95) $C(\gamma)$ in configuration 1\\
\L(.57*.95) $C(\gamma)$ in configuration 2\\
\L(.10*.71) $Q^{1}[\gamma,\nu]$\\
\L(.77*.25) $Y^{2}[\gamma,\nu_1,\nu_2]$\\
\L(.25*.86) $\nu$\\
\L(.29*.33) $\eta$\\
\L(.55*.84) $\nu_1$\\
\L(.80*.84) $\nu_2$\\
\L(.67*.64) $p$\\
\L(.25*.47) $p$\\
\L(.17*.34) $\delta'$\\
\L(.63*.36) $\delta'$\\
\L(.34*.28) $\delta''$\\
\L(.69*.34) $\delta''$\\
\L(.23*.09) $\gamma$\\
\L(.66*.04) $\gamma$\\
\L(.27*.14) $\gamma'$\\
\L(.58*.09) $\gamma'$\\
\L(.75*.09) $\gamma''$\\
\endSetLabels
\AffixLabels{%
\centerline{%
\includegraphics[height=6cm,width=12cm]{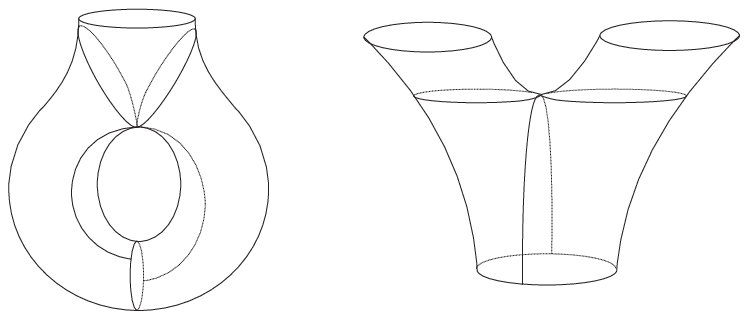}}}
\caption{$C(\gamma)$ in configuration 1 and 2}
\label{fig:collar}
\end{figure}

To prove the main theorems we will have to consider on many occasions the configuration in which the closure of the collar of a scg self-intersects. The closure of the collar of a scg $\gamma$, $\overline{C(\gamma)}$ self-intersects in a single point $p$. There exist two geodesic arcs of length $w$ emanating from $\gamma$ and perpendicular to $\gamma$ having the endpoint $p$ in common. These two arcs, $\delta'$ and $\delta''$,  form a smooth geodesic arc $\delta$. Two cases are possible - either $\delta$ arrives at $\gamma$ on opposite sites of $\gamma$ or it arrives on the same side (\textit{see Fig.~\ref{fig:collar}.}).

\begin{defi}  The closure of the collar of a scg $\gamma$, $\overline{C(\gamma)}$ self-intersects in a point $p$. We say that \textbf{$C(\gamma)$ is in configuration 1} if the shortest geodesic arcs $\delta'$ and $\delta''$ emanating from the intersection point $p$ and meeting $\gamma$ perpendicularly arrive at $\gamma$ on opposite sides. We say, that \textbf{$C(\gamma)$ is in configuration 2}, if they arrive on the same side of $\gamma$.
\end{defi}

For both configurations we have a corresponding Y-piece, a topological three-holed sphere, whose interior is isometrically embedded in $S$. If $C(\gamma)$ is in configuration 1, we cut open $S$ along $\gamma$. We call $S'$ the surface obtained in this way from $S$. Let $\gamma^1$ and $\gamma^2$ the two scg in $S'$ corresponding to $\gamma$ in $S$. Let $\nu$ be the shortest scg in the free homotopy class of $\gamma^1 \delta \gamma^2 \delta^{-1}$. Then $\gamma^1$,$\gamma^2$ and $\nu$ bound a three-holed sphere $Y^{1}$, whose interior lies in $S'$. As this decomposition occurs frequently, we will refer to it as $Y^{1}[\gamma,\nu ]$, \textit{the Y-piece for $\gamma$ from configuration 1}. If we close $Y^{1}[\gamma,\nu ]$ again at $\gamma$, we obtain a R.S. of signature $(1,1)$, $Q^{1}[\gamma,\nu ] \subset S$ (\textit{see Fig.~\ref{fig:collar}.}). Note, that in this case we obtain
\begin{equation}
\nu < 2\gamma + 2\delta = 2\gamma +4w,
\label{eq:Y1bound}
\end{equation}

as $\nu$ is in the free homotopy class of $\gamma^1 \delta \gamma^2 \delta^{-1}$. However, we can also calculate the exact value of $\nu$ by decomposing $Y^{1}[\gamma,\nu ]$ into two isometric hexagons, $H_1$ and $H_2$. Here we cut open $Y^{1}[\gamma,\nu ]$ along the shortest geodesic arcs connecting the boundary curves. In $H_1$ $\delta$ is the shortest geodesic arc connecting $\frac{\gamma^1}{2}$ and $\frac{\gamma^2}{2}$ and $\frac{\nu}{2}$ is the side opposite of $\delta$. From the geometry of right-angled hexagons (see \cite{b}, p. 454) we obtain

\[
 \cosh(\frac{\nu}{2}) = \sinh(\frac{\gamma}{2})^2  \cosh(\delta) - \cosh(\frac{\gamma}{2})^2
\]

As  $\cosh(x)^2 = \sinh(x)^2+ 1$ this is equal to

\begin{equation}
 \nu = 2\arccosh(\sinh(\frac{\gamma}{2})^{2}(\cosh(2w)-1) -1).
\label{eq:Y1bound_small}
\end{equation}

We note furthermore that there exists a geodesic arc $\gamma'$ in $\gamma$ connecting the two endpoints of $\delta$ on $\gamma$, whose length is restricted by $\gamma' \leq \frac{\gamma}{2}$. The shortest scg $\eta$ in the free homotopy class of $\gamma' \delta$ has length smaller than
\begin{equation}
\eta < \frac{\gamma}{2} + \delta = \frac{\gamma}{2} +2w .
\label{eq:Y1bound_two}
\end{equation}

If $C(\gamma)$ is in configuration 2, $\delta$ divides $\gamma$ in two parts, $\gamma'$ and $\gamma''$. Let $\nu_1$ and $\nu_2$ be the scg in the free homotopy class of $\gamma' \delta$ and $\gamma'' \delta$. The three scg $\gamma$, $\nu_1$ and $\nu_2$ then bound a Y-piece, we will refer to it as $Y^{2}[\gamma,\nu_1, \nu_2]$, \textit{the Y-piece for $\gamma$ from configuration 2} (\textit{see Fig.~\ref{fig:collar}.}). Note that $\nu_1 < \gamma' + \delta$ and $\nu_2 < \gamma'' + \delta$. Let WLOG $\gamma' \leq \gamma''$. As $\gamma = \gamma' \cup  \gamma''$, we have
\begin{eqnarray}
    \nu_1 < \gamma' + \delta \leq \frac{\gamma}{2}+\delta = \frac{\gamma}{2}+2w \text{ \ \ and \ \ }
    \nu_2 < \gamma'' + \delta < \gamma+2w.
\label{eq:Y2bound}
\end{eqnarray}

For small values of $\gamma$, we obtain a better bound for $\nu$ by decomposing $Y^{2}[\gamma,\nu_1, \nu_2]$ into two isometric hexagons, $H_1$ and $H_2$, by cutting it open along the shortest geodesic arcs connecting the boundary curves. Here $\frac{\delta}{2}$ is the unique geodesic arc in $H_1$ perpendicular to $\frac{\gamma}{2}$ and the geodesic arc between $\frac{\nu_1}{2}$ and $\frac{\nu_2}{2}$ and with endpoints on both arcs. $\frac{\delta}{2}$ divides $H_1$ into two pentagons, $P_1$ and $P_2$. Let $P_1$ be the pentagon, that contains $\frac{\gamma'}{2}$ as a boundary arc. From the geometry of right-angled pentagons (see \cite{b}, p. 454), we get
\begin{equation}
\cosh(\frac{\nu_1}{2})= \sinh(\frac{\gamma'}{2})\sinh(\frac{\delta}{2})   \text{ \ \ and \ \ }  \nu_1 \leq 2\arccosh(\sinh(\frac{\gamma}{4})\sinh(w)),
\label{eq:Y2bound_small}
\end{equation}

as $\sinh$ and $\arccosh$ are strictly increasing functions on $\mathbb{R}^{+}$.

With the help of this decomposition, the following lemma is proven in \cite{bs}, p. 40-42 :

\begin{lem}
Let $S$ be a compact R.S. of genus $g$ and $\gamma$ a scg in $S$. Let $C(\gamma)$ be the collar of $\gamma$ of width $w$. If $C(\gamma)$ is in configuration 1, let $\delta$ be the geodesic arc emanating from both sides of $\gamma$ and perpendicular to $\gamma$. $\delta$ divides $\gamma$ into two arcs. Let $\gamma'$ be the shorter of the two. Let furthermore $\eta$ be the scg in the free homotopy class of $\gamma' \delta$. If $\eta \geq \gamma$, then
\[
   w  \geq \max\left\{\arcsinh\left(\frac{1}{\sinh(\frac{\gamma}{2})}\right),\arccosh\left( \frac{\cosh(\frac{\gamma}{2})}{\cosh(\frac{\gamma}{4})}\right)\right\} \geq W'
\]
If $C(\gamma)$ is in configuration 2, let $Y^{2}[\gamma,\nu_1, \nu_2]$ be the Y-piece for $\gamma$ from configuration 2. If either $\nu_1$ or $\nu_2$ is bigger than $\gamma$, then
\[
    w \geq \arccosh(2)= W
\]
\label{thm:collar_bound}
\end{lem}

The lower bound for the width of the geodesic $\gamma$ depends on the constant $K$, where
\[
   K = 3.326.
\]
As a consequence of this lemma we have, that if $C(\gamma)$ is in configuration 1 and $\gamma < 3.326 = K$ then its width $w$ is bigger than $\arctanh(2/3)=W'$. If $C(\gamma)$ is in configuration 1 and $\gamma > K$ and bigger than $\arccosh(2)=W$. If $C(\gamma)$ is in configuration 2 its width $w$ is always bigger than $W$.\\

For the proof of \textbf{Theorem~\ref{thm:main}} we also need the following lemma from \cite{bs}, p. 38 :

\begin{lem}
Let $F$ be a compact Riemann surface of signature $(h,1)$, such that $1 \leq h $ and assume, that the boundary $\eta$ of $F$ has length $\eta < 2\log(8h-2)$  Then $F$ contains a nsscg $\alpha$ of length smaller than $2\log(8h-2)$ in its interior.
\label{thm:lem_bound1}
\end{lem}

A consequence of this lemma is, that every compact R.S. of genus $g$ contains a nsscg of length smaller than $2\log(4g-2)$ in its interior (see \cite{bs}, p. 38). With the help of this lemma, we prove the following :

\begin{lem}
Let $F$ be a compact Riemann surface of signature $(h,1)$ and assume, that the boundary $\eta$ of $F$ has length $\eta$. Then $F$ contains a nsscg $\alpha$ of length smaller than \\
$L=\max \left\{\frac{\eta}{2}+\log(8h-2),2\log(8h-2) \right\}$ in its interior.
\label{thm:lem_bound2}
\end{lem}

\textbf{proof of Lemma~\ref{thm:lem_bound2}} The collar of $\eta$, $C(\eta)$ of width $w$ in $F$ is in configuration 2. Let $Y^2[\eta,\nu_1,\nu_2]$ be the Y-piece for $\eta$ from configuration 2 and $\nu_1 \leq \nu_2$. We have by~\eqref{eq:Y2bound} \\
\[
\nu_1 \ < \frac{\eta}{2} + 2w
\]
We now show that either $F$ contains a nsscg $\alpha$ of length $2\log(8h-2)$ in its interior or that $2w < \log(8h-2)$, from which follows, that $\nu_1 < L$. If $\nu_1$ is non-separating, then we are done. If not, we cut open $F$ along $\nu_1$ into two parts. The part $F^1$, that does not contain $\eta$ has signature $(h^1,1)$, where $h^1 \leq h-1$ and its boundary is $\nu_1 \leq  L $. In this case we argue as before and divide $F^1$ again into two parts. As long as the shorter scg in the Y-piece for the boundary geodesic from configuration 2 is separating, we can successively cut off pieces $F^k$ from $F$. Let $(h^k,1)$ be the signature of $F^k$, where $h^k \leq h-k$. Repeating the argument for $\nu_1$ we obtain that the boundary geodesic of a $F^k$ has length smaller than $L$. This procedure ends at least, when $F^k$ is a Q-piece, a Riemann surface of signature (1,1). Then the decomposition of $F^k$ yields a nsscg $\alpha$ of length smaller than $L$ in the interior of $F^k \subset F$.\\
To conclude the proof, we have to show, that $2w < \log(8h-2)$.
Consider the surface $F'=F + F / \eta$, which is obtained by attaching the mirror image of $F$ along the boundary $\eta$. It has genus $2h$. As a consequence of \textbf{Lemma~\ref{thm:lem_bound1}} there exists a nsscg $\alpha$ of length smaller than $2\log(8h-2)$ in the interior of of $F'$. Note, that $\alpha \neq \eta$, as $\eta$ is separating in $F'$. If $\alpha \cap \eta = \emptyset$, then $\alpha$ is contained in $F$ and we are done. If $\alpha \cap \eta  \neq \emptyset$, then it has to traverse the collar of $\eta$, $C(\eta)$ in $F'$ at least twice and therefore $2\log(8h-2) > \alpha > 4w$, from which follows, that $2w < \log(8h-2)$. $\square$ \\

With the help of the previous lemmata we establish an upper bound for the second shortest scg on a compact R.S. in the following lemma. Other methods were applied in \cite{b}, p. 123 to obtain such an upper bound, however the one obtained here is lower.

\begin{lem}
Let $S$ be a compact Riemann surface of genus $g \geq 2$ and let $\gamma_1$ be the systole of $S$ and $\gamma_2$ be the second shortest scg on $S$. Then $\gamma_1 \leq 2\log(4g-2)$ and $\gamma_2 \leq 3\log(8g-7)$.
\label{thm:gamma2}
\end{lem}

\textbf{proof of Lemma~\ref{thm:gamma2}} By an area argument (see \cite{b}, p. 124) the length of the shortest scg, $\gamma_1$ of a compact Riemann surface $S$ of genus $g$ is bounded from above by $2\log(4g-2)$. If $\gamma_1$ is separating, we cut open $S$ along $\gamma_1$ into two parts $S^1$ and $S^2$. Let WLOG $S^1$ be the part, such that $S^1$ is of signature $(h,1)$, such that $h \leq \frac{g}{2}$. By \textbf{Lemma~\ref{thm:lem_bound1}} there exists a nsscg $\alpha$ of length smaller or equal to $2\log(4g-2)$ in the interior of $S^1$. In this case we have $\gamma_2 \leq 2\log(4g-2)$.\\
If $\gamma_1$ is non-separating, we have to take another approach. The collar of $\gamma_1$, $C(\gamma_1)$ intersects in the point $p_1$ and has width $w_1$. We furthermore know, that the interior of $C(\gamma_1)$ is isometrically embedded into $S$ and therefore its area can not exceed the area of $S$. Therefore
\[
2 \gamma_1 \cdot \sinh(w_1) = \area C(\gamma_1) < \area S=4\pi(g-1)
\]

Hence
\[
w_1 \leq \arcsinh(\frac{2\pi(g-1)}{\gamma_1})
\]
If $\frac{\pi}{2} \leq \gamma_1 \leq 2\log(4g-2)$, we obtain an upper bound for $w_1$, using that $\arcsinh(x) \leq \log(2x+1)$.
\[
 w_1 \leq \log(8(g-1)+1) < \log(8g-7)
\]
In this case, we can conclude, that there is a scg $\gamma_2 \neq \gamma_1$ in $S$ of length smaller than
\[
    \gamma_2 < \frac{\gamma_1}{2} + 2 w_1  < 3\log(8g-7).
\]
To see this, we apply either equation (~\ref{eq:Y1bound_two}) or equation (~\ref{eq:Y2bound}), depending on, whether the collar of $\gamma_1$ is in configuration 1 or 2, respectively.\\
If $\gamma_1 < \frac{\pi}{2}$, we have to consider again both possible configurations. If $\gamma_1$ is in configuration 1, we obtain by equation (~\ref{eq:Y1bound_small}), using the decomposition of the Y-piece from configuration 1, $Y^{1}[\gamma_1,\nu ]$ into hexagons and as $w_1 \leq \arcsinh(\frac{2\pi(g-1)}{\gamma_1})$ that
\[
 \nu \leq 2\arccosh((\sinh(\frac{\gamma_1}{2})^{2}(\cosh(2\arcsinh(\frac{2\pi(g-1)}{\gamma_1}))-1)) -1) \leq 4\log(8g-7).
\]
Here the upper bound of $4\log(8g-7)$ was determined using MAPLE. If we cut open $S$ along $\nu$, we obtain two pieces, one corresponding to $Y^{1}[\gamma_1,\nu ]$ and a second piece $S'$ of signature $(g-1,1)$. Applying \textbf{Lemma~\ref{thm:lem_bound2}} to $S'$, we conclude, that there exists a scg $\gamma_2$ in $S' \subset S$, whose length is bounded from above by $\frac{\nu}{2} + \log(8(g-1)-2) \leq 3\log(8g-7)$.\\
If $\gamma_1$ is in configuration 2, we obtain from the decomposition of the Y-piece from configuration 2, $Y^{2}[\gamma_1,\nu_1,\nu_2 ]$ into pentagons (equation (~\ref{eq:Y2bound_small})) and as $w_1 \leq \arcsinh(\frac{2\pi(g-1)}{\gamma_1})$ that
\[
\nu_1 \leq 2\arccosh(\sinh(\frac{\gamma}{4})\frac{2\pi(g-1)}{\gamma_1}) \leq 3\log(8g-7).
\]
Again the upper bound of $3\log(8g-7)$ was determined using MAPLE. $\square$ \\

A useful result for Riemann surfaces with boundary was obtained in \cite{ge} :

\begin{lem}
Let $S$ be a Riemann surface of signature $(g,n)$. Let $\gamma_1$ be the systole of $S$ and $l(\partial S)$ be the length of the boundary of $S$. Then $\gamma_1 \leq 4\log(4g+2n+3) + l(\partial S)$.
\label{thm:boundary}
\end{lem}

For a Q-piece, a R.S. of signature $(1,1)$ we have the following inequalities for a short canonical homology basis $(\ao,\at)$ by \cite{pa1}, p. 59-62 :

\begin{lem}
Let $Q$ be a Riemann surface of signature $(1,1)$ and $\gamma$ be the boundary geodesic of $Q$. There exists a canonical homology basis $(\ao,\at)$ , $\ao \leq \at$  of $Q$ satisfying the following inequalities :\\
\[
    \cosh(\frac{\ao}{2}) \leq \cosh(\frac{\gamma}{6})+\frac{1}{2}  \text{ \ \ and \ \ }
\]
\[
    \cosh(\frac{\at}{2}) \leq \sqrt{\frac{\cosh^2 (\frac{\gamma}{4})+\cosh^2(\frac{\ao}{2})-1}{2(\cosh(\frac{\ao}{2})-1)}}.
\]
\label{thm:Qpiece}
\end{lem}

The result is stated differently in \cite{pa1}. In \cite{pa1} $\ao$ is the shortest scg in the interior of $Q$ and $\at$ the shortest scg in $Q$, that intersects $\ao$. But due to this construction, both $\ao$ and $\at$ are non-separating. Furthermore $\at$ intersects $\ao$ only once due to its minimality. Hence $\ao$ and $\at$ have the required properties for a canonical homology basis.\\
Another lemma needed for the proof of the main theorem  concerns comparison surfaces and can be found in \cite{pa2}, p. 234 :
\begin{lem}
Let $S$ be a Riemann surface of signature $(g,n)$ with $n > 0$. Let $\beta_1,...,\beta_n$ be the boundary geodesics of $S$. For $(\epsilon_1,...,\epsilon_n) \in (\mathbb{R}^+)^n$ with at least one $\epsilon_i > 0$, there exists a comparison surface $S_c$ with boundary geodesics of length $\beta_{1}+\epsilon_1,...,\beta_{n}+\epsilon_n$ such that for each simple closed geodesic $\gamma_c$ in the interior of the comparison surface $S^c$, there exists a geodesic $\gamma$ in the interior of $S$, such that $\gamma < \gamma^c$.
\label{thm:comp_S}
\end{lem}

We finally state a consequence of the collar lemma stated in \cite{bs}, p. 106 :

\begin{lem}
Let $S$ be a Riemann surface of genus $g$ with $g \geq 2$. Let $\gamma$ be a simple closed geodesic in $S$. If $\eta$ is another scg, that does not intersect $\gamma$, then
\[
      \arcsinh(\frac{1}{\sinh(\frac{\gamma}{2})}) < \dist(\eta,\gamma)
\]
and if $w$ is the width of the collar of $\gamma$, $C(\gamma)$, then $w > \arcsinh(\frac{1}{\sinh(\frac{\gamma}{2})})$.
\label{thm:col_lem}
\end{lem}

\section{Main theorems}

\begin{figure}[h!]
\SetLabels
\L(.10*.95) \textbf{Case 1.} $\gamma_1 \neq \ao$\\
\L(.42*.92) $\gamma_1$\\
\L(.31*.88) $\ao$\\
\L(.51*.90) $\at$\\
\L(.10*.75) \textbf{Case 2.} $\gamma_1 = \ao$\\
\L(.10*.69) \textbf{Case 2.a.1.)}\\
\L(.20*.64) $\gamma_2$\\
\L(.29*.62) $\ao$\\
\L(.09*.60) $\at$\\
\L(.55*.69) \textbf{Case 2.a.2.)}\\
\L(.64*.64) $\gamma_2$\\
\L(.49*.58) $\ao$\\
\L(.62*.57) $\at$\\
\L(.10*.44) \textbf{Case 2.b.1.)} $\gamma_2 = \at$\\
\L(.29*.38) $\ao$\\
\L(.33*.33) $\at$\\
\L(.55*.44) \textbf{Case 2.b.2.)} $\gamma_2 \neq \at$\\
\L(.49*.33) $\ao$\\
\L(.64*.32) $\at$\\
\L(.58*.28) $\gamma_2$\\
\L(.10*.19) \textbf{Case 2.c.1.)}\\
\L(.27*.00) $\gamma_2$\\
\L(.27*.17) $\ao$\\
\L(.44*.10) $\at$\\
\L(.19*.08) $\beta_1$\\
\L(.35*.08) $\beta_2$\\
\L(.55*.19) \textbf{Case 2.c.2.)}\\
\L(.71*.00) $\gamma_2$\\
\L(.71*.17) $\ao$\\
\L(.75*.10) $\at$\\
\L(.62*.08) $\beta_1$\\
\L(.78*.08) $\beta_2$\\
\endSetLabels
\AffixLabels{%
\centerline{%
\includegraphics[height=12cm,width=14cm]{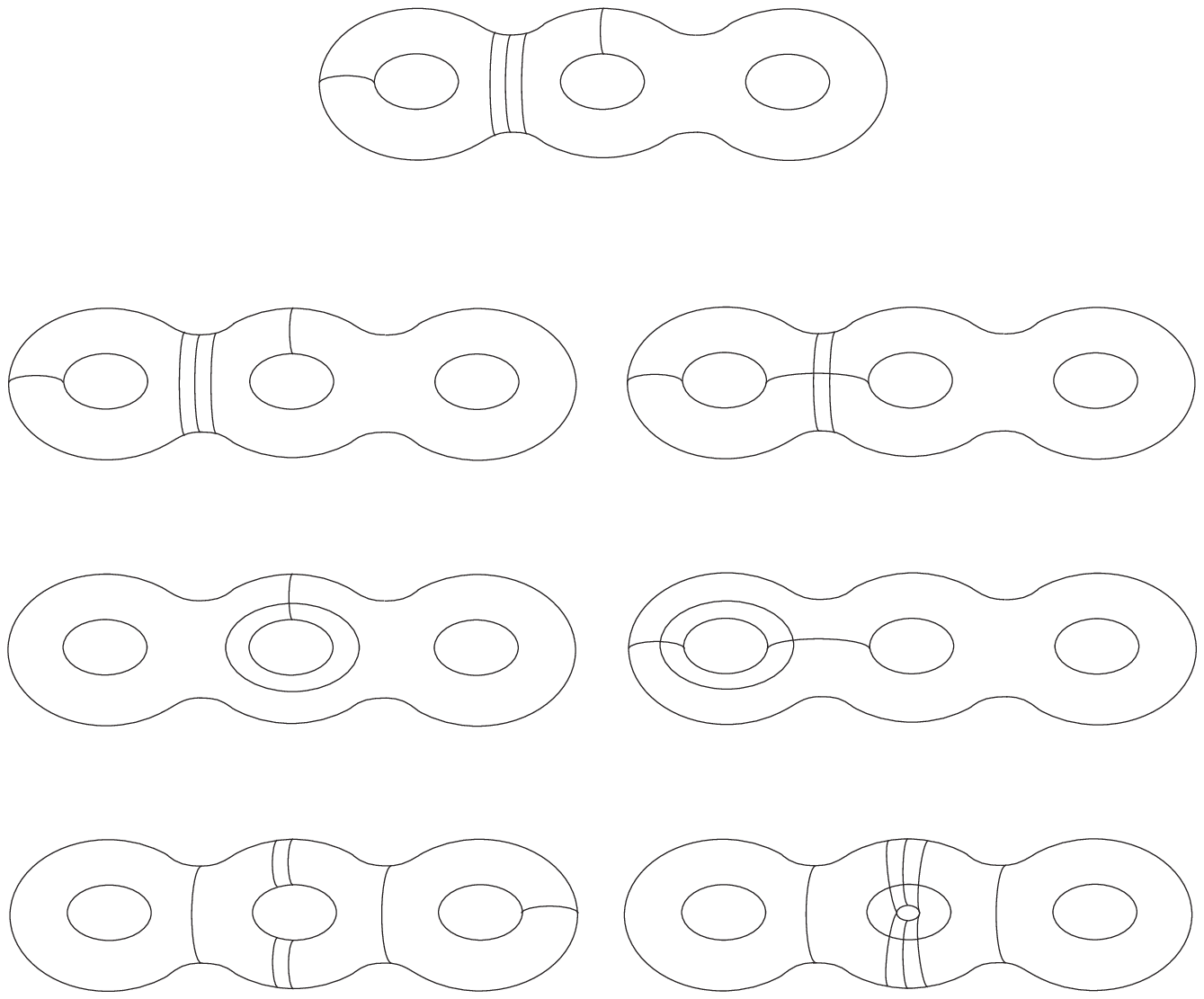}}}
\caption{Relative positions of $\ao$ and $\at$ in the different cases of the proof of Theorem~\ref{thm:main}}
\label{fig:cases}
\end{figure}

\textbf{proof of Theorem~\ref{thm:main}} It is well known, that two nsscg $\ao$ and $\at$ can be incorporated together into a canonical homology basis, if $\ao \cup \at$ does not separate $S$ into two parts and if $\ao$ and $\at$ have either exactly one or no intersection point. To prove  \textbf{Theorem~\ref{thm:main}} we have to show, that there exist two short nsscg, $\ao$ and $\at$, with these properties and whose collar width is bounded from below. Then we can obtain \textbf{Theorem~\ref{thm:main}} from equation (\ref{eq:mk2}). The proof of \textbf{Theorem~\ref{thm:main}} depends on whether the shortest scg $\gamma_1$ in $S$, the systole, is separating or non-separating. We will distinguish several cases. These cases are depicted in Fig.~\ref{fig:cases}.

\textbf{Case 1.} : \textit{The systole $\gamma_1$ of $S$ is a separating scg} \\

By  \textbf{Lemma~\ref{thm:gamma2}}, $\gamma_1$ has length smaller than $2 \log(4g-2)$. We cut open $S$ along $\gamma_1$, which yields two R.S., $S_1$ and $S_2$ of signature $(h_1,1)$ and $(h_2,1)$, respectively, such that $h_1 \leq h_2$. In both surfaces the half-collar of $\gamma_1$ is in configuration 2. By \textbf{Lemma~\ref{thm:collar_bound}} the width of a half-collar of $\gamma_1$ is bigger than $W$. We now show, that there exist two short nsscg, $\ao$ in $S_1$ and $\at$ in $S_2$, whose collars in $S$ have width of at least $W$. As $\ao \cup \at$ cannot divide $S$ into two parts and as $\ao$ and $\at$ do not intersect, they can be both together incorporated into a canonical homology basis of $S$. By \textbf{Lemma~\ref{thm:lem_bound1}} and \textbf{Lemma~\ref{thm:lem_bound2}}, we have that

\[
\ao < 2 \log(4g-2)   \text{ \ \ and \ \ } \at < 2 \log(8(g-1)-2)=2\log(8g-10)
\]

We now show, that each $\alpha_i, i\in \{1,2\}$ has a collar, whose width $w_i$ is bounded from below. Namely, if $\alpha_i < K$ then  $w_i > W'$ and if $\alpha_i > K$, then $w_i > W$.

\begin{figure}[h!]
\SetLabels
\L(.45*.73) $\delta^*$\\
\L(.47*.65) $s_1$\\
\L(.47*.37) $s_2$\\
\L(.56*.55) $p'$\\
\L(.47*.99) $\gamma_1'$\\
\endSetLabels
\AffixLabels{%
\centerline{%
\includegraphics[height=5cm,width=3cm]{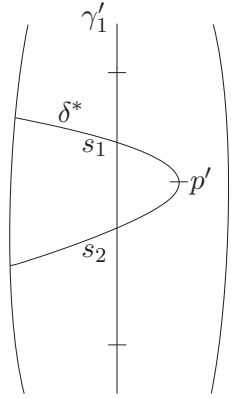}}}
\caption{Lift of $C(\gamma_1)$ in the universal covering}
\label{fig:colint}
\end{figure}

Consider WLOG the collar of $\ao$ in $S_1$. Its closure self-intersects in a point $p \in S_1$ or a geodesic arc $\delta$ of length smaller than $w_1$, emanating perpendicularly from $\ao$ meets the boundary of $S_1$ first. \\
In the first case, we apply \textbf{Lemma~\ref{thm:collar_bound}}.  We obtain, that if  $\ao < K$ then  $w_1 > W'$ and if $\ao > K$, then $w_2 > W$. In the second case, we show, that $\overline{C(\ao)}$ can not self-intersect in $C(\gamma_1) \cap S_2$. Therefore it self-intersects a point $p \in S_2 \backslash C(\gamma_1)$. In this case every geodesic arc $\delta'$ emanating perpendicularly from $\ao$ with endpoint $p$ has to traverse $C(\gamma_1) \cap S_2$ and hence has length bigger than $W$ in $S$.\\
To prove, that $\overline{C(\ao)} \subset S$ can not self-intersect in $S_1 \cup  C(\gamma_1)$, we lift $C(\gamma_1)$ into the universal covering. Here $\gamma_1$ lifts to $\gamma_1'$ and $\delta \cap C(\gamma_1)$ to $\delta^*$ (\textit{see Fig.~\ref{fig:colint}.}). The lift $\delta^*$ is a geodesic. Let $s_1$ and $s_2$ be the first intersection points of $\delta^*$ and $\gamma_1'$, the lift of $\gamma_1$ on opposite sides of $p'$, the lift of $p$. There exists an unique geodesic arc connecting $s_1$ and $s_2$, which is an arc in $\gamma_1'$, as $\gamma_1'$ is a geodesic. But $s_1$ and $s_2$ also lie on $\delta^*$, which implies, that $\delta^*$ is contained in $\gamma_1'$, a contradiction.\\

\textbf{Summary of Case 1:}
If the systole $\gamma_1$ of $S$ is a separating scg, then we can always find two short nsscg $\ao < 2\log(4g-2)$ and $\at < 2\log(8g-10)$ for a homology basis of $S$. Let $w_1$ and $w_2$ be the collar width of $\ao$ and $\at$, respectively. If  $\alpha_i < K$ then  $w_i > W'$ and if $\alpha_i > K$, then $w_i > W$, for $i \in \{1,2\}$. It follows from equation (~\ref{eq:mk1}) and (~\ref{eq:mk2}) and the subsequent remark that $m_1(J(S))^2$ and $m_2(J(S))^2$ satisfy the inequalities from \textbf{Theorem~\ref{thm:main}}.\\

\smallskip

\textbf{Case 2.} : \textit{The systole $\gamma_1$ of $S$ is a non-separating scg }\\

In this case we can find a homology basis of $S$, such that $\gamma_1=\ao$. As $\ao$ is the shortest nsscg, it follows from equation (~\ref{eq:mk1}) and (~\ref{eq:mk2}) that $m_1(J(S))^2$ satisfies the inequalities from \textbf{Theorem~\ref{thm:main}}.\\
To find a second short scg, that does not separate $S$ together with $\ao$, we have to consider several subcases. Let $\gamma_2$ be the second shortest scg on $S$. By \textbf{Lemma~\ref{thm:gamma2}} its length smaller than $3\log(8g-7)$. We will have to examine different cases, depending on whether $\gamma_2$ is separating, non-separating and non-separating with $\ao$ or non-separating but separating together with $\ao$.\\

\textbf{ Case 2.a)}\textit{ $\gamma_2$ is separating}\\

Note, that $\gamma_1$ and $\gamma_2$ can not intersect. It is easy to see, that otherwise we could find a scg in $S$, that is smaller than $\gamma_2$. We separate $S$ into two parts, $S_1$ and $S_2$ along $\gamma_2$. Let $S_1$ be the part, which contains $\ao$ and $S_2$ be the remaining part of signature $(h_2,1)$, such that $h_2 \leq g-1$. In this case $\gamma_2$ is smaller than $2\log(8g-10)$, due to the minimality of $\gamma_2$. Otherwise we would arrive again at a contradiction, if we apply \textbf{Lemma~\ref{thm:lem_bound2}} to $S_2$.
The collar of $\gamma_2$ is in configuration 2. Let $Y^2[\gamma_2,\nu_1,\nu_2]$ be the Y-piece for $\gamma_2$ from configuration 2. We have to distinguish two cases for the choice of $\at$, where the choice depends on $Y^2[\gamma_2,\nu_1,\nu_2]$.\\

\textbf{ Case 2.a.1.)}\textit{ $Y^2[\gamma_2,\nu_1,\nu_2] \neq Y^2[\gamma_2,\ao,\ao]$}\\

In this case let $\at$ be the shortest nsscg in $S_2$. As $\gamma_2 < 2\log(8g-10)$ it follows from \textbf{Lemma~\ref{thm:lem_bound2}} that $\at < 2\log(8g-10)$. $\ao$ and $\at$ can be incorporated together into a canonical homology basis. As $\ao$ does not occur twice in the boundary curves of $Y^2[\gamma_2,\nu_1,\nu_2]$ and as $\gamma_2$ is the second shortest scg in $S$, we conclude by \textbf{Lemma~\ref{thm:collar_bound}}, that the collar of $\gamma_2$ has width $w' \geq W$. We now determine a lower bound for the width of $C(\at)$, $w_2$. $\at$ is the shortest nsscg in the interior of $S_2$. Hence we can argue as in the case of the collar of $\at$ in \textbf{Case 1} to obtain a lower bound for the width of $C(\at)$, $w_2$. If $\at < K$ then  $w_2 > W'$ and if $\at > K$, then $w_2 > W$.\\

\smallskip

\textbf{ Case 2.a.2.)}\textit{ $Y^2[\gamma_2,\nu_1,\nu_2] = Y^2[\gamma_2,\ao,\ao]$}\\

If $\nu_1=\nu_2=\ao$, then the interior of $Y^2[\gamma_2,\nu_1,\nu_2]$ is embedded in the Q-piece $Q_1=S_1$, a R.S. of signature $(1,1)$. This case can not occur, if  $2.1 \leq \ao = \gamma_1 \leq \gamma_2$, because otherwise there would exist a scg $\at' \neq \ao$ in $Q_1$, that is smaller than $\gamma_2$ by \textbf{Lemma~\ref{thm:Qpiece}}.\\
In this case let $\beta$ be the shortest nsscg in $S_2$. As $\gamma_2 < 2\log(8g-10)$ it follows from \textbf{Lemma~\ref{thm:lem_bound2}} that $\beta < 2\log(8g-10)$. Let $\at$ be the shortest nsscg in $S$, that does not intersect $\ao$. We have $\at \leq \beta < 2\log(8g-10)$.\\
$\at$ has a collar $C(\at)$, whose width $w_2$ is bounded from below. To see this, we cut open $S$ along $\ao$ to obtain $S'$. Consider the collar of $\at$ in $S'$. Its closure self-intersects in a point $p \in S'$ or a geodesic arc $\delta$ of length smaller than $w_2$, emanating perpendicularly from $\at$ meets the boundary of $S'$ first. \\
By \textbf{Lemma~\ref{thm:col_lem}}
$\dist(\ao,\at) > \arcsinh(\frac{1}{\sinh(\frac{\ao}{2})}) > \arcsinh(\frac{1}{\sinh(\frac{2.1}{2})})$, as $\ao \leq 2.1$. It follows from the same arguments as in \textbf{Case 1.}, that $w_2$ has the lower bound
\[
w_2 > \min \left\{\arcsinh(\frac{1}{\sinh(\frac{2.1}{2})}),W' \right\} > 0.73.
\]
\\
\textbf{Summary of Case 2.a.) :} We can always find two short nsscg $\ao$ and $\at$ for a homology basis of $S$,whose lengths satisfy the same upper bounds as in as in \textbf{Case 1.} and whose collar width is bounded from below, such that $m_1(J(S))^2$ and $m_2(J(S))^2$ satisfy the inequalities from \textbf{Theorem~\ref{thm:main}}.\\

\smallskip

\textbf{ Case 2.b)} \textit{$\gamma_2$ is non-separating and non-separating with $\gamma_1=\ao$}\\

In this case we have to distinguish two cases,  $\at=\gamma_2$ and  $\at \neq \gamma_2$.\\

\textbf{ Case 2.b.1.)} $\at=\gamma_2$\\

We have $\at = \gamma_2 < 3\log(8g-7)$. Note, that $\at$ can not intersect $\ao$ more than once, as otherwise there would exist a scg, which is shorter than $\at$. We now determine a lower bound for the width of the collar of $\at$, $C(\at)$.\\
If $C(\at)$ is in configuration 2, let $Y^{2}[\at,\nu_1, \nu_2]$ be the Y-piece for $\at$ from configuration 2. If $\nu_1$ and $\nu_2$ are both smaller than $\at$, then both must be $\ao$. If $Y^{2}[\at,\nu_1, \nu_2]=Y^{2}[\at,\ao,\ao]$, then $Y^{2}[\at,\ao,\ao]$ is embedded in $S$ as a Q-piece with boundary $\at$ and $\at$ would be separating, a contradiction. Hence we conclude by \textbf{Lemma~\ref{thm:collar_bound}} the width of $C(\at)$ is bigger than $W$.\\

\begin{figure}[h!]
\SetLabels
\L(.29*.99) $\at'$\\
\L(.69*.99) $\at^*$\\
\L(.37*.73) $s'$\\
\L(.38*.61) $\theta$\\
\L(.34*.56) $r_2$\\
\L(.56*.54) $\delta_{2}'$\\
\L(.53*.36) $\lambda'$\\
\L(.66*.25) $s^*$\\
\endSetLabels
\AffixLabels{%
\centerline{%
\includegraphics[height=5cm,width=6cm]{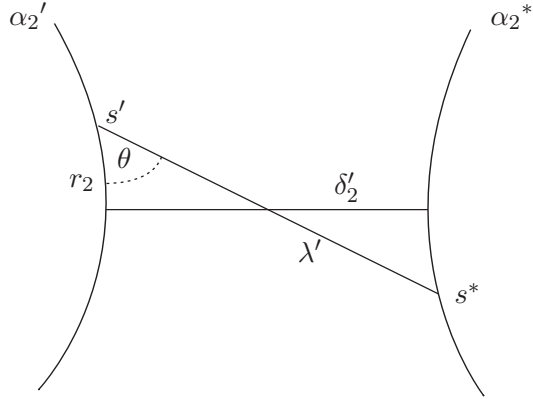}}}
\caption{Two lifts of $\at$ in the universal covering}
\label{fig:collar_atone}
\end{figure}

If $C(\at)$ is in configuration 1, $\overline{C(\at)}$ self-intersects in a single point $p$. There exist two geodesic arcs of length $w_2$ emanating from $\at$ and perpendicular to $\at$ having the endpoint $p$ in common. These two arcs form a smooth geodesic arc $\delta_2$. We lift $\at$ and $\delta_2$ in the universal covering. Here $\at$ lifts to $\at'$ and $\at^*$ and $\delta_2$ to $\delta_2'$. In the covering there exist two points, $s' \in \at'$ and $s^* \in \at^*$, on opposite sites of $\delta_2'$ and at the same distance $r_2 \leq \frac{\at}{4}$ from $\delta_2'$, such that $s'$ and $s^*$ are mapped to the same point $s \in \at$ by the covering map. By drawing the geodesic $\lambda'$ from $s'$ to $s^*$, we obtain two isometric right-angled geodesic triangles. (\textit{see Fig.~\ref{fig:collar_atone}.})  \\

We have to consider two subcases, $\lambda' \neq \ao$ and $\lambda'=\ao$ and $\ao >1.1$. In the case $\lambda'=\ao \leq 1.1$, we will switch to \textbf{Case 2.b.2.)} \\

\textit { $\lambda' \neq \ao$ }\\

If $\lambda' \neq \ao$, then we can again argue as in \textbf{Case 1}. We obtain, that if  $\at < K$ then  $w_2 > W'$ and if $\at > K$, then $w_2 > W$.\\

\textit{ $\lambda'=\ao$ and $\ao >1.1$ }\\

To intersect $\ao$, $\at$ has to traverse the collar of $\ao$. We can use this fact to derive a lower bound for the width of the collar of $\at$, $C(\at)$.

\begin{figure}[h!]
\SetLabels
\L(.47*.99) $\ao'$\\
\L(.47*.68) $q_2$\\
\L(.62*.68) $u_2$\\
\L(.55*.55) $\at''$\\
\L(.47*.50) $s''$\\
\L(.48*.43) $\theta$\\
\L(.51*.42) $r_1$\\
\L(.36*.33) $u_1$\\
\L(.51*.33) $q_1$\\
\endSetLabels
\AffixLabels{%
\centerline{%
\includegraphics[height=6cm,width=4.5cm]{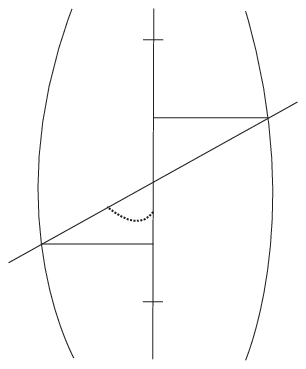}}}
\caption{Lift of $C(\ao)$ in the universal covering}
\label{fig:collar_ao}
\end{figure}

Lift $C(\ao)$ in the universal covering and let $C'(\ao)$ be the lift of $C(\ao)$ (\textit{see Fig.~\ref{fig:collar_ao}}). $\at$ traverses the collar $C(\ao)$ of width $w_1$. It lifts to $\at''$ and in the lift it enters $C'(\ao)$ at a point $u_1$ and leaves at a point $u_2$. Consider the geodesic arcs emanating from $u_1$ and $u_2$ respectively and meeting the lift of $\ao$, $\ao'$ perpendicularly. Their length is $w_1$. Let $q_1$ and $q_2$ be the endpoints of these geodesic arcs on $\ao'$. $\at''$ intersects $\ao'$ in the midpoint $s''$ of the geodesic arc between $q_1$ and $q_2$ under angle $\theta$. Here $s''$ is a lift of $s \in Q_1$. Set $r_1= \dist(q_1,s'')=\dist(q_2,s'')$. Then $r_1$ is smaller or equal to $\frac{\ao}{4}$, as $\at$ is the third shortest scg in $Q_1$ and otherwise there exists another point $u_2'$, such that $u_2$ and $u_2'$ map to the same point on $Q_1$ under the universal covering map, such that $\dist(u_1,u_2') < \dist(u_1,u_2)$, a contradiction to the fact, that $\at$ is minimal.
Consider the right-angled triangle with vertices $u_1$, $q_1$ and $s''$. From the geometry of hyperbolic triangles (see \cite{b}, p. 454) we have for $\theta$ :

\[
     \sin(\theta)=\frac{\sinh(w_1)}{\sinh(\dist(u_1,s''))}  \text{ \ \ and \ \ } \cosh(\dist(u_1,s''))= \cosh(w_1)\cdot \cosh(r_1)
\]
from which follows, as $\sinh^2(x)=\cosh^2(x)-1$, that
\[
    \sin(\theta)=\frac{\sinh(w_1)}{\sqrt{\cosh^2(r_1) \cdot \cosh^2(w_1)-1}}
\]
The point $s''$ corresponds to $s'$ in the other lift of $\at$ (\textit{see Fig.~\ref{fig:collar_atone}.}) and the angle $\theta$ to the interior angle of the right-angled geodesic triangle at the vertex $s'$. From the geometry of this triangle we get :

\[
     \sin(\theta)=\frac{\sinh(w_2)}{\sinh(\frac{\ao}{2})}
\]
and therefore, as $\sinh(w_2)$ is decreasing with increasing $r_1 \leq \frac{\ao}{4}$
\begin{equation}
\sinh(w_2) \geq \frac{\sinh(w_1) \cdot \sinh(\frac{\ao}{2})}{{\sqrt{\cosh^2(\frac{\ao}{4}) \cdot \cosh^2(w_1)-1}}}
\label{eq:wtwo_a}
\end{equation}

Note that the left hand side in \eqref{eq:wtwo_a} is increasing with increasing $w_1$ and increasing $\ao$.
As the width of $C(\ao)$, $w_1$ is bigger than $W'$, we get a lower bound for $w_2$, if we set $w_1 = W'$. In this case we obtain from equation~\eqref{eq:wtwo_a}

\begin{equation}
 w_2 \geq {w_2}^Q = \arcsinh\left( \frac{\frac{2\sqrt{5}}{5} \cdot \sinh(\frac{\ao}{2})}{{\sqrt{\frac{9}{5}\cosh^2(\frac{\ao}{4})-1}}} \right)
\label{eq:qwtwo}
\end{equation}

In this case we obtain, with $\ao > 1.1$ :
\[
 m_2(J(S))^2 < \frac{3\log(8g-7)}{\pi-2 \cdot \arcsin(\frac{1}{\cosh({w_2}^Q)})} \leq 3.1\log(8g-7) .
\]

\textbf{Summary of Case 2.b.1.) :} We can always find two short nsscg $\ao = \gamma_1 < 2\log(4g-2)$ and $\at = \gamma_2 <3\log(8g-7)$ for a homology basis of $S$. Their collar width is bounded from below, such that $m_1(J(S))^2$ and $m_2(J(S))^2$ satisfy the inequalities from \textbf{Theorem~\ref{thm:main}}.\\

\smallskip

\textbf{ Case 2.b.2.)} $\at \neq \gamma_2$\\

This case treats the remaining case for $C(\gamma_2)$ in configuration 1, and the geodesic $\lambda'$ (see Fig.~\ref{fig:collar_atone}.) is $\ao$, with  $\ao  \leq 1.1$.\\
In this case we cut open $S$ along $\ao$ to obtain the surface $S'$ of signature $(g-1,2)$. Let $\ao'$ and $\ao''$ be the boundary.
In this case we let $\at$ be the shortest nsscg in $S'$, that does not intersect $\ao$. We first show the following claim.

\begin{cla}
The shortest nsscg $\at \subset S'$ has length smaller than $2\log(24g-23)+2.2$ .
\label{thm:scut_nsscg}
\end{cla}

\textbf{proof of \textbf{Claim~\ref{thm:scut_nsscg}}}  We first show, that there exists a scg of length smaller than $2 \log(24g-23)+2.2$ in the interior of $S'$. It is sufficient to proof this statement for the case $\ao =1.1$. It follows from \textbf{Lemma~\ref{thm:comp_S}}, that this is also true for $\ao < 1.1$. If there exists a scg of length smaller than $2 \log(24g-23)+2.2$ in $S'$ and this geodesic is non-separating, we are done. If it is separating, we apply \textbf{Lemma~\ref{thm:lem_bound2}} and conclude, that there exists a nsscg in $S'$, that is smaller than $\log(24g-23)+1.1 + \log(8g-10) < 2\log(24g-23)+2.2$, which proves the claim.\\
Let $\ao' = 1.1.$. The closure of the half-collar $\overline{C(\ao')} \subset S'$ self-intersects in a point in $S'$ or a geodesic arc emanating perpendicularly from $\ao'$, of length smaller than $w'$ meets $\ao''$ perpendicularly in a point $p_1$. We examine two cases, which depend on how $\overline{C(\ao')}$ intersects itself. \\

 \textit{i) The closure of the half-collar of $\ao'$ intersects $\ao''$ in $p_1$ before self-intersecting in $S'$\\}

A geodesic arc $\sigma \subset S'$  meets $\ao'$ and $\ao''$ perpendicularly on both endpoints where $p_1$ is the endpoint on $\ao''$. The distance set $Z_r(\ao') \subset S'$ is defined by
\[
Z_r(\ao') = \{ x \in S' \mid \dist(x,\ao') < r \}.
\]
As long as $r$ is small enough, such that $Z_r(\gamma) \subset C(\ao')$, we have from hyperbolic geometry that $\area(Z_r(\ao'))=\ao'  \cdot \sinh(r)$. Consider $Z_{\sigma}(\ao')$. It is embedded into $S'$ and therefore its area can not exceed the area of $S'=S$, which is smaller than $4\pi(g-1)$. Therefore
\[
   \ao'  \cdot \sinh(\sigma) =  \area Z_{\sigma}(\ao')   < \area(S) = 4\pi(g-1)
\]

As $\ao' =1.1 $ and as $\arcsinh(x) \leq \log(2x+1)$, we obtain an upper bound for $\sigma$.
\[
    \sinh(\sigma) \leq \frac{4\pi(g-1)}{1.1}   \Rightarrow  \sigma \leq \log(24g-23)
\]
Hence we conclude, that the shortest scg $\beta_1 \subset S'$ in the free homotopy class of $\ao' \sigma \ao'' \sigma^ {-1} $ has length smaller than $2\log(24g-23)+2.2$.\\

\textit{ii) The closure of the collar of $\ao$ self-intersects in $p_1 \in S'$}\\

A geodesic arc $\sigma$ passes through $p_1$ and meets $\ao$ perpendicularly on both endpoints. Let $Y^2[\ao,\nu_1,\nu_2]$ be the Y-piece for $\ao$ from configuration 2.\\
As $ \ao =1.1$, we conclude by the same area argument as in case \textit{i)}, that $\sigma <  2\log(24g-23)$. From equation ~\eqref{eq:Y2bound} it follows, that both $\nu_1$ and $\nu_2$ are smaller than
\[
\ao + \sigma \leq 2\log(24g-23)+1.1 .
\]
At least one of them is not $\ao''$. Hence there exists a scg of length smaller than $2\log(24g-23)+2.2$ in $S'$.
Hence we have proven the claim.  $\square$ \\

Let $w_2$ be the width of the collar of $\at$. In this case we conclude as in \textbf{ Case 2.a.2)}, that $\dist(\ao,\at) > \arcsinh(\frac{1}{\sinh(\frac{\ao}{2})}) > \arcsinh(\frac{1}{\sinh(\frac{1.1}{2})})$, as $\ao \leq 1.1$. It follows from the same arguments as in \textbf{Case 2.a.2.)}, that $w_2$ has the lower bound
\[
w_2 > \min \left\{\arcsinh(\frac{1}{\sinh(\frac{1.1}{2})}),W' \right\} > W'.
\]

\textbf{Summary of Case 2.b.2.) :} We have that $\at < 2\log(24g-23)+2.2$ and $w_2 >W'$. We obtain, that
\[
     m_2(J(S))^2 < \frac{2\log(24g-23)+2.2}{\pi-2\arcsin(\frac{1}{\cosh(W')})} \leq 3.1 \log(8g-7)
\]

\smallskip

\textbf{ Case 2.c)} \textit{$\gamma_2$ is non-separating, but separating with $\gamma_1 = \ao$}\\

By \textbf{Lemma~\ref{thm:gamma2}} we know, that the length of $\gamma_2$ is bounded by $3 \log(8g-7)$. It is easy to see, that $\gamma_2$ can not intersect $\ao$. It can not intersect $\ao$ more than once, due to the minimality of the two geodesics and it can not intersect $\ao$ once due to the fact, that it is separating with $\ao$. As $\gamma_2$ is separating with $\ao$, we conclude by \textbf{Lemma~\ref{thm:collar_bound}}, that its collar width is bounded from below. If  $\gamma_2 < K$ then  the width of its collar is bigger than $W'$ and if $\gamma_2 > K$, then the width of its collar is bigger than $W$.\\
We cut open $S$ along $\gamma_2$ and $\ao$. The two geodesics divide $S$ into $S_1$ and $S_2$. We first show, the following claim :

\begin{cla}
The shortest nsscg $\at^i \subset S_i$ has length smaller than $4.5 \log(8g-7)$ for $i \in \{1,2\}$.
\label{thm:s1_nsscg}
\end{cla}

The proof is similar to the proof of \textbf{Claim~\ref{thm:scut_nsscg}}.\\
\textbf{proof of \textbf{Claim~\ref{thm:s1_nsscg}}} Consider WLOG $S_1$. We proof the claim for the cases $\ao < \pi$ and $\ao \geq \pi$ :\\

\textit{a) $\ao \geq \pi$ \\}

$\overline{C(\ao)}$ self-intersects in a point in $S_1$ or a geodesic arc emanating perpendicularly from $\ao$, of length smaller than $w_1$ meets $\gamma_2$ perpendicularly in a point $p_1$. We examine two cases, which depend on how $\overline{C(\ao)}$ intersects itself. \\

 \textit{i) The closure of the collar of $\ao$ intersects $\gamma_2$ in $p_1$ before self-intersecting in $S_1$\\}

A geodesic arc $\sigma \subset S_1$  meets $\ao$ and $\gamma_2$ perpendicularly on both endpoints where $p_1$ is the endpoint on $\gamma_2$. We now define for a scg $\gamma$ in $S$ and an $r>0$, the distance set of distance $r$ of $\gamma$, $Z_r(\gamma)$ by
\[
Z_r(\gamma) = \{ x \in S \mid \dist(x,\gamma) < r \}.
\]
As long as $r$ is small enough, such that $Z_r(\gamma) \subset C(\gamma)$, we have from hyperbolic geometry that $\area(Z_r(\gamma))=2\gamma  \cdot \sinh(r)$. Consider $Z_{\sigma}(\ao) \cap S_1$. It is embedded into $S_1$ and therefore its area can not exceed the area of $S_1$, which is smaller than $4\pi((g-2)-1)$. Therefore
\[
   \ao  \cdot \sinh(\sigma) =  \area Z_{\sigma}(\ao) \cap S_1  < \area S_1=4\pi(g-3).
\]

As $ \pi \leq  \ao $ and as $\arcsinh(x) \leq \log(2x+1)$, we obtain an upper bound for $\sigma$.
\[
    \sinh(\sigma) \leq \frac{4\pi(g-3)}{\pi}   \Rightarrow  \sigma \leq \log(8(g-3)+1) < \log(8g-7)
\]
Hence we conclude, that the shortest scg $\beta_1$ in the free homotopy class of $\ao \sigma \gamma_2 \sigma^ {-1} $ has length smaller than $7 \log(8g-7)$. It is a separating scg. Applying \textbf{Lemma~\ref{thm:lem_bound2}} we conclude, that there exists a nsscg of length smaller than $ 4.5  \log(8g-7)$ in $S_1$. Note that, using the hexagon decomposition (see \cite{bs}, p. 454) of the Y-piece with boundary geodesics $\beta_1$, $\ao$ and $\gamma_2$, we can obtain the exact value of the length of $\beta_1$, which will be useful later for small values of $\ao$. It is

\begin{equation}
\cosh(\frac{\beta_1}{2}) = \sinh(\frac{\ao}{2})\sinh(\frac{\gamma_2}{2})\cosh(\sigma)-\cosh(\frac{\ao}{2})\cosh(\frac{\gamma_2}{2}).
\label{eq:beta1}
\end{equation}

\textit{ii) The closure of the collar of $\ao$ self-intersects in $p_1 \in S_1$\\}

A geodesic arc $\sigma$ passes through $p_1$ and meets $\ao$ perpendicularly on both endpoints. Let $Y^2[\ao,\nu_1,\nu_2]$ be the Y-piece for $\ao$ from configuration 2.\\
if $ \ao \geq  \pi $, we conclude by the same area argument as in case \textit{i)}, that $\sigma <  \log(8g-7)$. From equation ~\eqref{eq:Y2bound} it follows, that both $\nu_1$ and $\nu_2$ are smaller than
\[
\ao + \sigma \leq 3 \log(8g-7).
\]
At least one of them is not $\gamma_2$. Therefore, if this scg is non-separating, we are done. If it is separating, we cut off the part of $S_1$, that contains $\ao$ and conclude by \textbf{Lemma~\ref{thm:lem_bound2}}, that this part contains a nsscg of length smaller than $3 \log(8g-7) < 4.5  \log(8g-7)$. This settles the claim in the case $\ao \geq \pi$.\\

\textit{b) $\ao < \pi$ \\}

If  $\ao < \pi$, we use the fact, that there exists a comparison surface $S_{1}^{c}$ for $S_1$, as described in \textbf{Lemma~\ref{thm:comp_S}}, such that one boundary geodesic has length $\pi$ and the other has length $\gamma_2$ and conclude, that it contains a scg of length smaller than $4.5 \log(8g-7)$ in its interior. Therefore there exists a scg of length smaller than $4.5\log(8g-7)$ in $S_1$, by \textbf{Lemma~\ref{thm:comp_S}}. If this geodesic is separating, we apply again \textbf{Lemma~\ref{thm:lem_bound2}} and conclude, that there exists a nsscg in $S_1$, that is smaller than $4.5\log(8g-7)$. $\square$ \\

In total, we obtain, that the shortest nsscg $\at^1$ in $S_1$ and the shortest nsssg $\at^2$ in $S_2$ are both smaller than $ 4.5 \log(8g-7)$. Both can be incorporated with $\ao$ into a canonical homology basis. Consider the sets $Z_{W'}(\ao)$ and $Z_{W'}(\gamma_2)$, with $W'=\arctanh(2/3)$. We now choose a nsscg $\at$, that is non-separating with $\ao$. The choice depends on how  $Z_{W'}(\ao)$ and $Z_{W'}(\gamma_2)$ intersect. We distinguish two cases.\\

\textbf{Case 2.c.1.)} \textit{$Z_{W'}(\ao) \cap Z_{W'}(\gamma_2) \cap S_1 = \emptyset$ or $Z_{W'}(\ao) \cap Z_{W'}(\gamma_2) \cap S_2 = \emptyset$}\\

If $Z_{W'}(\ao) \cap Z_{W'}(\gamma_2) \cap S_1 = \emptyset$, then we choose $\at = \at^2 \subset S_2$ and if $Z_{W'}(\ao) \cap Z_{W'}(\gamma_2) \cap S_2 = \emptyset$ we choose $\at = \at^1 \subset S_1$. Consider WLOG the first case. We show, that the collar of $\at^2$ has width bigger than $W'$. If the closure of the collar of $\at^2$ self-intersects in $S_1$, it has to traverse either $S_1 \cap Z_{W'}(\ao)$ or $S_1 \cap Z_{W'}(\gamma_2)$ and hence its width is bigger than $W'$. This follows from the same arguments as in \textbf{Case 1}. If $\overline{C(\at^2)}$ self-intersects in $S_2$, we conclude by \textbf{Lemma~\ref{thm:collar_bound}}, that $\at^2$ has a collar of with bigger than $W'$.\\

\textbf{Summary of Case 2.c.1.) :} $\at$ is the shortest nsscg in either $S_1$ or $S_2$, its length is restricted by
$\at < 4.5 \log(8g-7)$, the width of its collar is bigger than $W'$. It follows from equation (~\ref{eq:mk1}) and (~\ref{eq:mk2}), that
\[
m_2(J(S))^2 < 3.1\log(8g-7).
\]

\smallskip

\textbf{Case 2.c.2.)} \textit{$Z_{W'}(\ao)$ and $Z_{W'}(\gamma_2)$ intersect both in $S_1$ and $S_2$} \\

\begin{figure}[h!]
\SetLabels
\L(.30*.90) $Y'$\\
\L(.70*.90) $Y''$\\
\L(.48*.90) $\ao$\\
\L(.16*.50) $\beta_1$\\
\L(.30*.50) $\delta'$\\
\L(.53*.50) $\at$\\
\L(.62*.50) $\delta''$\\
\L(.82*.50) $\beta_2$\\
\L(.48*.10) $\gamma_2$\\
\L(.25*.00) $F$\\
\endSetLabels
\AffixLabels{%
\centerline{%
\includegraphics[height=5cm,width=10cm]{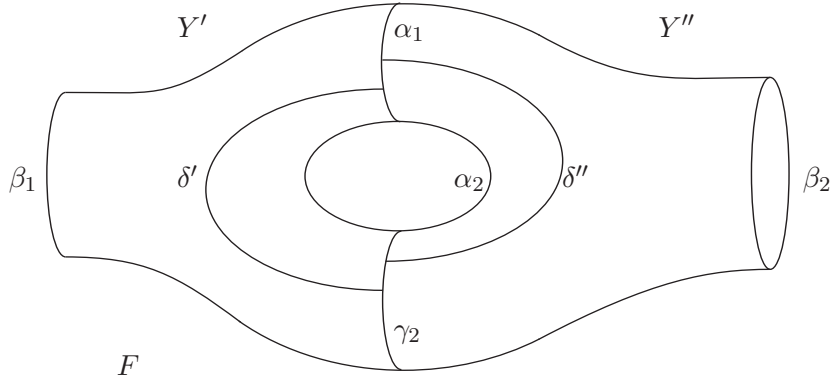}}}
\caption{The Riemann surface $F$ of signature (1,2)}
\label{fig:twosep}
\end{figure}

If $Z_{W'}(\ao)$ and $Z_{W'}(\gamma_2)$ intersect both in $S_1$ and $S_2$ we have to argue in a different way. We choose another small nsscg to be $\at$. We may assume, that $\ao \geq 1.5$ and $\gamma_2 \geq 2.1$. Otherwise it would follow from equation (~\ref{eq:beta1}), with $\sigma=2W'$, that $\beta_1 < \gamma_2$, a contradiction.\\
We now choose $\at$. Let $\delta'$ and $\delta''$ be the shortest geodesic arcs in $S_1$ and $S_2$, respectively, connecting $\ao$ and $\gamma_2$. Their length is bounded from above by $2W'$ as $Z_{W'}(\ao)$ and $Z_{W'}(\gamma_2)$ intersect. The endpoints of $\delta'$ and $\delta''$ divide each $\ao$ and $\gamma_2$ into two geodesic arcs. Let $\ao^*$ and $\gamma^*_2$ be the shorter of these arcs. We define $\at$ to be the shortest scg in the free homotopy class of $\delta' \ao^* \delta'' \gamma^*_2$. It intersects $\ao$ only once. The length of $\at$ is restricted by
\[
     \at < \frac{\ao}{2} +\frac{\gamma_2}{2} + 4W' < 2.5 \log(8g-7) + 4W'.
\]
Let $Y' \subset S_1$ be the Y-piece with boundaries $\ao$, $\gamma_2$ and the shortest scg $\beta_1$ in the free homotopy class of  $\ao \delta' \gamma_2 \delta'^ {-1} $. Let $Y'' \subset S_2$ be the Y-piece constructed in the same way in $S_2$, having as third boundary $\beta_2$.
The union $Y' \cup Y''$ in $S$ is embedded as a Riemann surface $F$ of signature (1,2) \textit{(see Fig. ~\ref{fig:twosep}.)}.

\begin{figure}[h!]
\SetLabels
\L(.75*.95) $\gamma_{2}'$\\
\L(.36*.86) $s^{**}$\\
\L(.47*.81) $s'$\\
\L(.45*.57) $\delta'$\\
\L(.52*.56) $\eta'$\\
\L(.59*.54) $q''$\\
\L(.74*.49) $\alpha_{1}'$\\
\L(.47*.46) $q'$\\
\L(.56*.29) $\delta''$\\
\L(.59*.18) $s''$\\
\L(.71*.24) $s^{*}$\\
\L(.29*.05) $\gamma_{2}''$\\
\endSetLabels
\AffixLabels{%
\centerline{%
\includegraphics[height=6cm,width=8cm]{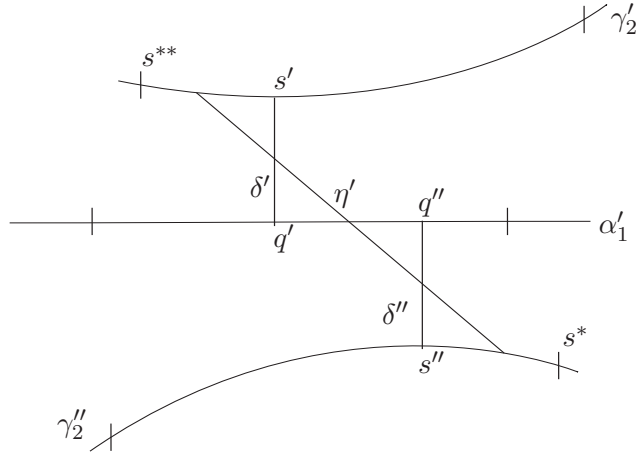}}}
\caption{Lifts of $\at$ and $\gamma_2$ in the universal covering}
\label{fig:twosep_det}
\end{figure}

However, if the length of $\gamma_2$ is small, the upper bound for $\at$ given above is not sufficient to establish an appropriate lower bound for the collar of $\at$. Therefore we will establish a better upper bound for the length of $\at$.\\
Lift $\ao$ to $\ao'$ and equally $\delta'$ and $\delta''$ into the universal covering \textit{(see Fig.~\ref{fig:twosep_det}.)}. By abuse of notation we will denote the lift of these two arcs by the same letter. To $\delta'$ and $\delta''$ attach the adjacent lifts of $\gamma_2$, $\gamma_2'$ and $\gamma_2''$ on opposite sides of $\ao'$. Let $q'$ be the endpoint of $\delta'$ on $\ao'$ and $q''$ be the endpoint of $\delta''$ on $\ao'$, such that $\dist(q',q'') \leq \frac{\ao}{2}$. Let furthermore be $s'$ be the endpoint of $\delta'$ on $\gamma_2'$ and $s''$ be the endpoint of $\delta''$ on $\gamma_2''$. Let  $s^{*}$ be the point on $\gamma_2''$, that maps to the same point on $\gamma_2$ under the covering map as $s'$, such that $\dist(s'',s^{*}) \leq \frac{\gamma_2}{2}$. Let equally be $s^{**}$ be the point on $\gamma_2'$, that maps to the same point on $\gamma_2$ under the covering map as $s''$, such that $\dist(s^{**},s')=\dist(s'',s^{*})$. Let $\eta'$ be the geodesic arc connecting the midpoint of $s'$ and $s^{**}$ on $\gamma'_2$ and midpoint of $s''$ and $s^{*}$ on $\gamma''_2$.\\

The image of $\eta'$ under the covering map, $\eta$  forms a closed geodesic arc on $S$. As $\eta$ is in the same free homotopy class as $\at$, its length provides an upper bound for the length of $\at$. In Fig.~\ref{fig:twosep_det}., the points $s^{**}$ and $s^{*}$ lie on opposite sides of $\delta'$ and $\delta''$. We will derive an upper bound for this case. In any other case the length of $\eta'$ is either shorter or the situation is a mirror image of the depicted one. It is clear, that $\eta'$ is maximal, if $\dist(s^{**},s')=\dist(s'',s^{*})=\frac{\gamma_2}{2}$ and $\dist(q',q'')=\frac{\ao}{2}=\frac{\gamma_2}{2}$. Therefore it is sufficient to discuss this case.\\
In this case we obtain from the geometry of hyperbolic triangles :
\begin{equation}
   \cosh(\frac{\at}{4}) < \cosh(\frac{\eta'}{4})=\cosh(\frac{\gamma_2}{4}) \cdot \cosh(\frac{\delta'}{2}) \leq  \cosh(\frac{\gamma_2}{4})\cosh(W')
\label{eq:eta_bound}
\end{equation}
as $\delta' \leq 2W'$. We now determine a lower bound for the width of the collar of $\at$. We have to distinguish several subcases :\\

\textbf{Case 2.c.2.a)} \textit{The collar of $\at$ is in configuration 1}  \\

$\overline{C(\at)}$ has width $w_2$. It self-intersects in a point $p$. Lift $\at$ into the hyperbolic plane as described \textbf{Case 2.a)} (\textit{see Fig.~\ref{fig:collar_atone}.}) We have to discuss two cases, $\lambda'=\ao$ and $\lambda' \neq \ao$. \\

\textit{$\lambda'=\ao$ }\\

This case was discussed in \textbf{Case 2.b)}. We may assume, that $\ao \geq 1.5$ from equation (~\ref{eq:beta1}). Hence we can apply equation (~\ref{eq:wtwo_a}) with $\alpha_1=1.5$ and $w_1=W'$ and obtain
\[
      w_2 > 0.66
\]

\textit{$\lambda' \neq \ao$ }\\

If $\lambda' \geq \at$ we can apply \textbf{Lemma~\ref{thm:collar_bound}} and conclude, that $w_2 \geq W'$. If $\lambda' < \at$, we conclude from equation(~\ref{eq:eta_bound}) and as $\gamma_2$ is the second shortest scg in $S$, that
\begin{equation}
  \gamma_2   \leq  \lambda' < \at < 4\arccosh(\cosh(\frac{\gamma_2}{4})\cosh(W')).
\label{eq:trlam}
\end{equation}
From the geometry of right-angled hyperbolic triangles (see \cite{b}, p. 454), we obtain from Fig.~\ref{fig:collar_atone}., that
\[
 \cosh(\frac{\at}{4})  \cosh(w_2)  \geq  \cosh(r_2) \cosh(w_2) = \cosh(\frac{\lambda'}{2}).
\]
Using the upper bound for $\at$ and the lower bound for $\lambda'$ from equation (~\ref{eq:trlam}) in this inequality we obtain :
\begin{equation}
     \cosh(w_2)  \geq    \frac{\cosh(\frac{\gamma_2}{2})}{\cosh(\frac{\gamma_2}{4})\cosh(W')}.
\label{eq:cf_one_spec}
\end{equation}

\textbf{Case 2.c.2.b)} \textit{The collar of $\at$ is in configuration 2}  \\

Let $Y^2[\at,\nu_1,\nu_2]$ be the Y-piece for $\at$ in configuration 2. $\overline{C(\at)}$ self-intersects in the point $p_2$, such that $\dist(p_2,\at)=w_2=\frac{\delta_2}{2}$. The geodesic arc $\delta_2$ emanating perpendicularly from $\at$ passes through this point and its endpoints divide $\at$ into two parts, $\at'$ and $\at''$. The common perpendiculars of the boundary geodesics of $Y^2[\at,\nu_1,\nu_2]$ separate the Y-piece into two isometric hexagons and $\delta_2$ decomposes these hexagons into pentagons. By the pentagon formula (see \cite{b}, p. 454) we have
\[
    \sinh(\frac{\delta_2}{2}) \sinh(\frac{\at'}{2})=\cosh(\frac{\nu_1}{2}) \text{ \ \ and \ \ }    \sinh(\frac{\delta_2}{2}) \sinh(\frac{\at''}{2})=\cosh(\frac{\nu_2}{2})
\]
None of these boundary geodesics can be $\ao$, as $\ao$ intersects $\at$. If either $\nu_1$ or $\nu_2$ is bigger, than $\at$, we obtain from \textbf{Lemma~\ref{thm:collar_bound}}, that $w_2 > W$. If not, then both must be bigger than $\gamma_2$. Additionally $\frac{\at'}{2}+ \frac{\at''}{2}=\frac{\at}{2}$. Therefore either $\frac{\at'}{2}$ or $\frac{\at''}{2}$ is bigger than $\frac{\at'}{4}$. Let WLOG $\at'$ be the bigger one. We obtain from equation ~\eqref{eq:eta_bound}:
\[
         \sinh(w_2)  \sinh(\arccosh(\cosh(\frac{\gamma_2}{4})\cosh(W'))) \geq \sinh(\frac{\delta_2}{2}) \sinh(\frac{\at'}{2})=\cosh(\frac{\nu_1}{2}) \geq \cosh(\frac{\gamma_2}{2}).
\]
or equally , as $\sinh(x)= \sqrt{\cosh^2(x)-1}$ for $ x\geq 0$ :
\begin{equation*}
    \sinh(w_2) \geq \frac{\cosh(\frac{\gamma_2}{2})}{\sqrt{\cosh^2(\frac{\gamma_2}{4})\cosh^2(W')-1}}.
\end{equation*}

It follows from this equation, that in \textbf{Case 2.c.2.b)}
\[
      w_2 > 0.96 .
\]

\textbf{Summary of Case 2.c.2) :}
$\at$ to be the shortest scg in the free homotopy class of $\delta' \ao^* \delta'' \gamma^*_2$ \textit{(see Fig.\ref{fig:twosep})}. Its length is restricted by
\[
\at < 4\arccosh(\cosh(\frac{\gamma_2}{4})\cosh(W')) .
\]
From the discussion of the subcases \textbf{Case 2.c.2.a)} we conclude, that the width of the collar $w_2$ is bounded from below by
\[
   w_2 \geq \min\left\{0.66,\arccosh\left( \frac{\cosh(\frac{\gamma_2}{2})}{\cosh(\frac{\gamma_2}{4})\cosh(W')}\right)\right\}
\]
As a consequence of equation (~\ref{eq:beta1}),we have that $2.1 \leq \gamma_2$. With the help of this lower bound it follows from the above equation, that $w_2$ is bounded from below. As $\at$ is bounded from above, it follows from  equation (~\ref{eq:mk2}) and (~\ref{eq:mk1}), that $m_2(J(S))^2$ is bounded from above. A refined analysis shows, that
\[
m_2(J(S))^2 < 3.1\log(8g-7).
\]

\bigskip

This proves that \textbf{Theorem~\ref{thm:main}} is valid. $\square$\\

\textbf{proof of Corollary~\ref{thm:cor_main}} The proof is very similar to the discussion of \textbf{Case 1.} of \textbf{Theorem~\ref{thm:main}}.

Let $\eta_i \leq t$ be one of the simple closed geodesics that divide $S$. By \textbf{Lemma~\ref{thm:col_lem}} the width of a half-collar of $\eta_i$ is bigger than $\arcsinh(\frac{1}{\sinh(\frac{t}{2})})$ on both sides of $\eta_i$. It follows also from the collar theorem, that any other scg in $S$ has a distance greater than $\arcsinh(\frac{1}{\sinh(\frac{t}{2})})$ from $\eta_i$. Let $S_i$ be a surface of genus $(g_i,n_i)$, $g_i >0$ from the decomposition of $S$. Let WLOG $\eta_1,...\eta_{n_i}$ be its boundary geodesics.  We first prove, that the shortest nsscg, $\alpha_i$ in $S_i$ is smaller than $(n_i+1)\max \{4\log(4g_i+2n_i+3),t\}$. Then we show, that it has a collar in $S$ whose width is bounded from below. By \textbf{Lemma~\ref{thm:boundary}}, there exists a scg $\gamma_i^1$ in $S_i$ of length  $\gamma_i^1  \leq 4\log(4g_i+2n_i+3) + n_i t$.\\
We have, that
\[
4\log(4g_i+2n_i+3) + n_i t \leq  (n_i+1)\max \{4\log(4g_i+2n_i+3),t\}.
\]
Hence, if $\gamma_i^1$ is non separating, we are done. If $\gamma_i^1$ is separating, we cut open $S_i$ along $\gamma_i^1$. $S_i$ decomposes into two surfaces, such that one of these two, $S_i^2$ has signature $(g'_i,n'_i)$, with $g'_i >0$ and $n'_i \leq n_i-1$. The length of its  boundary is smaller than $4\log(4g_i+2n_i+3) + (n_i -1)t$. We can again apply \textbf{Lemma~\ref{thm:boundary}} to this surface to obtain an upper bound for the length of a scg in $S_i$. Repeating this process iteratively we obtain, that there exists a nsscg in $S_i$, whose length is smaller than $(n_i+1)\max \{4\log(4g_i+2n_i+3),t\}$.\\
Each $\alpha_i, i\in \{1,...,m\}$ has a collar, whose width $w_i$ is bounded from below. Namely, if $\alpha_i < K$ then
$w_i > \min \left\{\arcsinh(\frac{1}{\sinh(\frac{t}{2})}),W' \right\}$ and if $\alpha_i > K$, then
$w_i > \min\left\{\arcsinh(\frac{1}{\sinh(\frac{t}{2})}),W \right\}$.
This follows from the same line of argumentation as in \textbf{Case 1.} of \textbf{Theorem~\ref{thm:main}}.
The $(\alpha_i)_{i=1,...,m}$ can be together incorporated into a canonical homology basis of $S$. From the bounds on the length and the width of the collars of the geodesics follows the bound on the norm of the lattice vectors of the Jacobian of $S$.
In total we obtain \textbf{Corollary~\ref{thm:cor_main}}. $\square$\\

\textbf{proof of Theorem~\ref{thm:hyper}} For the proof of \textbf{Theorem~\ref{thm:hyper}} we first give a suitable definition of a hyperelliptic surface.
\begin{defi}
Let $S$ be a compact Riemann surface of genus $g$. An \textit{involution} is an isometry $\phi : S \rightarrow S, \phi \neq id$, such that $\phi^2=id$. The surface $S$ is \textit{hyperelliptic}, if it has an involution, that has exactly
$2g+2$ fixed points. These fixed points are called the \textit{Weierstrass points (WPs)}.
\end{defi}

\begin{figure}[h!]
\SetLabels
\L(.52*.25) $p_1$\\
\L(.52*.76) $p_1^*$\\
\L(.66*.76) $p_2^*$\\
\L(.66*.25) $p_2$\\
\L(.73*.25) $p_3$\\
\L(.81*.25) $p_4$\\
\L(.17*.25) $p_5$\\
\L(.25*.25) $p_6$\\
\L(.32*.25) $p_7$\\
\L(.46*.25) $p_8$\\
\L(.46*.90) $\ao$\\
\L(.46*.64) $\at$\\
\L(.53*.12) $\gamma$\\
\L(.56*.74) $\gamma'$\\
\L(.52*.90) $\gamma_1'$\\
\L(.52*.64) $\gamma_2'$\\
\L(.49*.46) $S$\\
\L(.48*.34) $S \backslash \phi$\\
\endSetLabels
\AffixLabels{%
\centerline{%
\includegraphics[height=5cm,width=10cm]{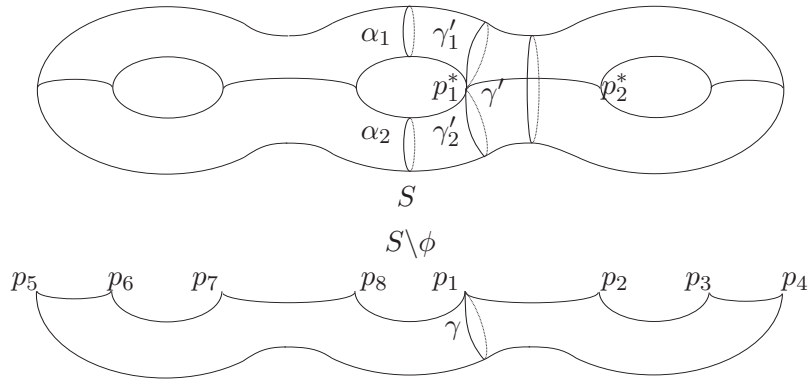}}}
\caption{A hyperelliptic surface $S$ and the quotient surface $S \backslash \phi$}
\label{fig:hyp}
\end{figure}

It is well known, that the above definition is equivalent to the usual one. Let $S$ be a hyperelliptic surface of genus $g$ with involution $\phi$. We will show, that the shortest nsscg $\ao$ of $S$ is bounded by a constant, independent of the genus. It was shown in \cite{bs}, that the width of the collar of the shortest nsscg on a R.S. $S$ is bounded from below. It then follows from equation (~\ref{eq:mk2}) and (~\ref{eq:mk1}), that \textbf{Theorem~\ref{thm:hyper}} holds.\\
Consider the quotient surface $S \backslash \phi$. This surface is a topological sphere with $2g+2$ cones of angle $\pi$, whose vertices $\{p_i\}_{i=1..2g+2}$ are the images of the WPs $\{p_{i}^*\}_{i=1..2g+2}$ under the projection \textit{(see Fig.~\ref{fig:hyp}.)}.\\

Let $B_r(p_i)$ be a disk of radius $r$ around a vertex of a cone. As long as these disks are small enough, they are embedded in $S\backslash \phi$. In this case the area of a disk of radius $r$ around a vertex of a cone $p_i$, $B_r(p_i)$ is half the area of a disk of radius $r$ in the hyperbolic plane, $\area ( B_r(p_i)) = \pi(\cosh(r)-1)$. Now expand all disks around the cone points until either a disk self-intersects or two different disks intersect for the first time at radius $R$. In this limit case we still obtain:
\[
  (2g+2)\pi(\cosh(R)-1) = \area\left(\bigcup\limits_{i = 1}^{2g + 2} {B_R (p_i )} \right)  < \area(S\backslash \phi) = 2\pi(g-1).
\]
As $\frac{g-1}{g+1} < 1$ we conclude, that $R < \arccosh(2)$.\\
When the radii of the disks reach $R$ and two different disks intersect the geodesic arc that forms lifts to a simple closed geodesic in $S$. When a disk self-intersects at radius $R$, the geodesic arc that forms lifts to a figure $8$ geodesic in $S$. This figure $8$ geodesic consists of two loops. The scg in the free homotopy class of such a loop is smaller than the loop itself. Hence there exists a scg of length smaller than $4R$ in $S$. It follows, for the systole $\gamma_1$ in $S$, that $ \gamma_1 < 4 \arccosh(2) = 5.2678...$.\\
By a refinement of this area estimate Bavard obtains a better upper bound in \cite{ba}, which is
\begin{equation}
 \gamma_1 < 4\arccosh \left( {\left(
2\sin \left( {\tfrac{{\pi (g + 1)}}
{{12g}}} \right)\right)}^{-1} \right) < 2\log(3+2\sqrt{3}+2\sqrt{5+3\sqrt{3}})=5.1067...
\label{eq:bav}
\end{equation}
We now show, that this upper bound is equally valid for the shortest non-separating scg in $S$. Consider the case, where two different disks intersect at radius $R$. In this case a geodesic arc of length smaller than $2R$ connects WLOG $p_1$ and $p_2$. It is easy to see, that it lifts to a scg $\ao$ of length $4R$ in the double covering $S$ \textit{(see Fig.~\ref{fig:hyp}.)}. $\ao$ passes the two WPs $p_{1}^*$ and $p_{2}^*$. It is a well-known fact, that such a geodesic is non-separating.
Consider now the case, where WLOG $B_R(p_1)$ self-intersects. The geodesic arc $\gamma$, that passes $p_1$ and the intersection point, lifts to a figure $8$ geodesic $\gamma'$ in $S$ \textit{(see Fig.~\ref{fig:hyp}.)}. Let $\gamma_1'$ and $\gamma_2'$ be the two different lifts of $\gamma$ in $S$ with intersection point $p_1^*$. Let $\ao$ and $\at$ be the scg in the free homotopy class of $\gamma_1'$ and $\gamma_2'$, respectively. The length of both is bounded from above by $2R$. It is easy to convince oneself, that the situation depicted Fig.~\ref{fig:hyp}. is the correct one and that both $\ao$ or $\at$ are non-separating.\\
In any case there exists a nsscg $\ao$ in $S$, whose length is smaller than the constant from equation (~\ref{eq:bav}). Hence we obtain an upper bound for $m_1(J(S))^2$ :
\[
   m_1(J(S))^2 <  \frac{3\log(3+2\sqrt{3}+2\sqrt{5+3\sqrt{3}})}{\pi} \leq 2.4382...
\]

This proves \textbf{Theorem~\ref{thm:hyper}}. $\square$

\end{document}